\documentclass[a4paper,11pt]{amsart}

\usepackage{epsfig}
\usepackage{color}
\usepackage{amssymb}

\newtheorem{theorem}{Theorem}
\newtheorem{corollary}[theorem]{Corollary}
\newtheorem{definition}{Definition}
\newtheorem{lemma}[theorem]{Lemma}
\newtheorem{proposition}[theorem]{Proposition}

\newtheorem*{proposition*}{Proposition}

\newcommand{\NN}{{\mathbb N}}

\newcommand{\RR}{{\mathbb R}}
\newcommand{\EU}{S}
\newcommand{\vv}{{\mathbf v}}
\newcommand{\ww}{{\mathbf w}}

\newcommand{\FF}{{\mathcal F}}
\newcommand{\GG}{{\mathcal G}_\tau}

\newcommand{\CC}{{\mathcal C}}

\newcommand{\ee}{{\mathrm e}}
\newcommand{\seta}{\longrightarrow}
\newcommand{\dpt}{\displaystyle}

\newcounter{lixo}

\DeclareMathOperator{\tr}{tr}

\DeclareMathOperator{\Fix}{Fix}

\title[Bifurcation diagrams of a forced cycle]{Transitions of bifurcation diagrams  of a forced heteroclinic cycle
\\ \today}

\author[I.S. Labouriau]{Isabel S. Labouriau}
\address{I.S. Labouriau\\
Centro de Matem\'atica
da Universidade do Porto \\
Rua do Campo Alegre,
687, 4169-007 Porto, Portugal }
\email{ islabour@fc.up.pt}

\author[A.A.P. Rodrigues]{Alexandre A. P. Rodrigues}
\address{A.A.P. Rodrigues\\ Center of Mathematics, Sciences Faculty, University of Porto \\ 
Rua do Campo Alegre 687, 4169--007 Porto, Portugal\\ 
and Lisbon School of Economics and Management \\
Rua do Quelhas 6, 1200-781 Lisboa, Portugal}
\email{alexandre.rodrigues@fc.up.pt, arodrigues@iseg.ulisboa.pt}

\thanks{
The first author has been partially supported by Centro de Matem\'atica da Universidade do Porto CMUP, member of LASI, which is financed by national  funds through  FCT/MCTES --- Funda\c{c}\~ao para a Ci\^encia e a Tecnologia, I.P. (Portugal)  under the projects with references  UIDB/00144/2020 and UIDP/00144/2020.}
\thanks{
The second author has been partially supported by the Project CEMAPRE/REM --- UIDB/05069/2020 financed
by FCT/MCTES (Portugal) through national funds.
}

\begin{document}

\begin{abstract}
A family of periodic perturbations of an attracting robust heteroclinic cycle defined on the two-sphere is studied by reducing the analysis to that of a one-parameter family of maps on a  circle.
The set of zeros of the family forms a bifurcation diagram on the cylinder.
The  different bifurcation diagrams and the transitions between them are obtained as the strength of attraction of the cycle and the amplitude of the periodic perturbation vary.
We determine a threshold in the cycle's attraction strength above which frequency locked periodic solutions with arbitrarily long periods bifurcate from the cycle as the period of the perturbation decreases.
Below this threshold further transitions are found giving rise to a frequency locked invariant torus and to a frequency locked suspended horseshoe, arising from heteroclinic tangencies in the family of maps.
 \end{abstract}

    \maketitle
  
  \bigbreak

\bigbreak
\textbf{2010 --- AMS Subject Classifications} 
{Primary: 37C60;   Secondary: 34C37, 34D20, 37C27, 39A23, 34C28}

\bigbreak
%\tableofcontents

\section{Introduction}\label{sec:intro}

Recent studies in several areas have emphasised ways in which \emph{heteroclinic cycles and networks} may be responsible for \emph{intermittent and non-regular dynamics} in nonlinear systems. They can be seen as the skeleton for the understanding of complicated dynamics. For example in \cite{Morrison22} the authors construct \emph{chaotic heteroclinic networks} to reproduce experimental data associated to intermittent behaviour of \emph{C. elegans},  a well studied model organism with simple anatomy, simple nervous system, few observed behaviours, and the relative ease of performing experimental measurements of   neural activity.

  A heteroclinic cycle in an autonomous continuous differential equation consists of saddle-type invariant sets and heteroclinic trajectories connecting them in a cyclic way. 
Such structures are \emph{structurally unstable} in generic dynamical systems without symmetry or other constraints. On the other hand, for equivariant vector fields the existence of invariant subspaces can force the existence of such connecting trajectories; heteroclinic cycles might then become \emph{robust} in the sense that they persist under small symmetry-preserving perturbations \cite{GH88, KM1, KM2, RL2014}.

Examples of robust heteroclinic cycles connecting equilibria have been discussed in many contexts. Perhaps the most important motivation for their study is connected with systems with heteroclinic cycles representing  the concept of \emph{winnerless competition} which has been widely discussed as being a more biologically-relevant paradigm than the alternative  \emph{winner-takes-all} assumption in many scenarios especially in Evolutionary Game Theory   \cite{Hofbauer98,  ML75, Rabinovich06}.

To date there has been   little systematic investigation of the effects of \emph{time-periodic perturbations}   and their \emph{bifurcations}, despite being natural for the modelling of many biological effects \cite{AHL2001, Mohapatra2015, WY2003}. 

Extending the results of \cite{LR18}, in this paper we study  the dynamics and bifurcations  in the neighbourhood of  a heteroclinic cycle contained in the asymptotically stable heteroclinic network explored in  \cite{ACL06} perturbed with a time-periodic forcing   with amplitude  $\gamma$.
  We reduce the problem to the analysis of the geometry of the set of zeros of a one-parameter family of functions defined on a circle---the \emph{bifurcation diagram}. 
We obtain a complete qualitative description of all possible persistent diagrams for different amplitudes of the forcing and different saddle-values of the network (a measure of its strength of attraction) as well as the transitions between different types. 
 This justifies the expression \emph{transitions of bifurcation diagrams} in the title.
We also find  possible integer multiples of the forcing frequency $\omega$ for which there is a  solution of the problem with frequency $n\omega$ for some $n\in\NN$ --- the \emph{frequency-locked} solutions.

Chaotic heteroclinic networks are good phenomenological models of switching dynamics \cite{PJ88, Rodrigues2013}. Given a homoclinic loop or heteroclinic cycle it has been known  that transversal intersections of stable and unstable manifolds give rise to complicated behaviour. The qualitative picture is clarified by Smale, who showed that this complicated behaviour included the existence of horseshoes.
 In the present article we identify 
 a discrete-time Bogdanov--Takens bifurcation  \cite{BRS} which can be seen as an  organising center to locate such horseshoes and strange attractors in the neighbourhood of
 the ghost of  the cycle that is broken for $\gamma>0$. 
  The study of \emph{singularities} can then be seen as a ``good route'' to locate chaos.

Although we use a specific example as model, our results can
be easily extended to a generic periodic perturbation of an autonomous differential equation with an attracting and clean heteroclinic  cycle in a network involving a finite number of saddles, provided the map associated to the first return is the same.  

\subsection*{Guide to the article}
The differential equation being perturbed is given in Section~\ref{sec:objectStudy}, where we present an overview of known results on its dynamics, as well as the form of perturbation.
The dynamics is codified by a discrete-time map on the cylinder and its expression is also provided in this section.

In Section~\ref{sec:bifCylinder} the analysis of the cylinder map is further reduced to the bifurcation of zeros of a family of maps of the circle parametrized by a time-lag parameter $\tau$.
Bifurcation diagrams with the distinguished parameter $\tau$ are  described 
 with a discussion of fold points
for different values of the strength $\delta$ of attraction of the heteroclinic cycle and of the amplitude $\gamma$ of the perturbation. The transitions between bifurcation diagrams as $\delta$ and $\gamma$ vary are also discussed.
The golden number $\delta=(1+\sqrt{5})/2=\Phi$ is obtained as a threshold for the diagrams. 

The stability of the bifurcating fixed points of the cylinder map around fold points is  studied in Section~\ref{sec:stability}.
Further dynamical transitions are also discussed and  a codimension two bifurcation is found for the weakly attracting case $\delta<\Phi$.
From this bifurcation it follows that there are parameter values  where an invariant circle is created
 and other values where horseshoe dynamics  arises  from tangencies of a saddle's invariant manifolds.

In Section~\ref{sec:freq} the results are interpreted back in terms of the original differential equations
 showing the occurrence 
 of frequency locking, invariant tori, periodic solutions with very long periods and chaotic regimes.
Section \ref{s: discussion} finishes the article with a discussion. 
We have endeavoured to draw illustrative bifurcation diagrams to make the paper easily readable.

\section{The object of study}\label{sec:objectStudy}
We are concerned with the 
 following set of the ordinary differential equations with a periodic forcing:
\begin{equation}
\label{general}
\left\{ 
\begin{array}{l}
\dot x = x(1-r^2)-\alpha x z +\beta xz^2 +\gamma(1-x)\sin(2\omega t)\\
\dot y =  y(1-r^2) + \alpha y z + \beta y z^2 \\
\dot z = z(1-r^2)-\alpha(y^2-x^2)-\beta z (x^2+y^2) 
\end{array}
\right.
\end{equation}
where 
$$
r^2= x^2+y^2+z^2, \qquad \beta<0<\alpha, \qquad  |\beta|<\alpha, \qquad \gamma, \omega\in \RR^+_0 .
$$
 The perturbing term $\gamma(1-x)\sin(2\omega t)$
 with amplitude controlled by $\gamma>0$
 appears only in the first coordinate
 because  it simplifies the computations and  allows comparison with previous work by other authors \cite{AHL2001, ACL06,DT3, Rabinovich06, TD1}.

\subsection{The unperturbed system ($\gamma=0$)}
\label{gamma=0}
 We denote by $F(x,y,z)$ the right hand side of the unperturbed system with $\gamma=0$ in \eqref{general}. 
This equation was  studied in \cite{ACL06}, we recall some of its properties. 

The unit sphere $\EU^2$ is  flow-invariant for \eqref{general} and attracts all trajectories except the origin which is a repelling equilibrium.  
The vector field $F$ has the following two symmetries of order 2
$$
\kappa_1(x,y,z)=(-x,y,z)
\qquad\mbox{and}\qquad
\kappa_2(x,y,z)=(x,-y,z) .
$$
 Only the symmetry $\kappa_2$ remains for  $\gamma>0$. 
The points  $\vv= (0,0,1)$  and $\ww= (0,0,-1)$ are equilibria.
From the symmetries  it follows that the planes $x=0$ and $y=0$ are flow-invariant 
and  meet $\EU^2$ in two flow-invariant circles connecting the equilibria $(0,0,\pm 1)$. 
Since  $\beta<0<\alpha$ and $\beta^2<8\alpha^2$, then these two equilibria are saddles and there is a pair of heteroclinic trajectories going from each equilibrium to the other. 
The derivative of $F$ at $p \in \{\vv, \ww\}$ has  expanding and contracting eigenvalues $E_p$ and $C_p$ given by:
$$
E_\vv= E_\ww= \alpha+\beta>0 \qquad \text{and} \qquad C_\vv= C_\ww= \beta-\alpha<0.
$$

The symmetry implies that the  invariant manifolds of $\vv$ and $\ww$  in $\EU^2$ coincide giving rise to a heteroclinic network $\Sigma_0$.  It consists of four cycles.
 In the restriction to each of the invariant planes $\Fix(\kappa_j)$, $j=1,2$,  the equilibria $\vv$ and $\ww$ have a saddle-sink connection, so this  network  is persistent under perturbations that preserve the symmetry and in this sense it is robust. An illustrative image of the network $\Sigma_0$ is given in \cite{RL2014}.

The constant $\delta^2$, called the saddle-value of $\Sigma_0$, measures the strength of attraction of the cycle in the absence of the periodic perturbation $\gamma(1-x)\sin(2\omega t)$. It  is given by
$$
\delta^2= \frac{C_\vv}{E_\vv} \frac{C_\ww}{E_\ww}= \frac{(\beta-\alpha)^2}{(\alpha+\beta)^2}>1.
$$
The Krupa--Melbourne criterion \cite{KM1, KM2} implies that $\Sigma_0$ 
 is asymptotically stable and, in particular, there are no periodic solutions near $\Sigma_0$.
 The case $\delta^2   = 1$ and $\gamma=0$ corresponds to a \emph{resonant bifurcation} of the robust heteroclinic cycle.
In the present article we assume  $\delta>1$.
 
 \subsection{The forced equation}\label{subsec:forced}
 The perturbed equation  
 $$
 \dot X=F(X)+(\gamma(1-x)\sin(2\omega t),0,0), \quad X=(x,y,z)\in \RR^3
 $$
 may be rewritten as an autonomous equation $(\dot X,\dot s)=\FF(X,s)$ in $\RR^3\times \EU^1$ 
 given by
 \begin{equation}\label{eq:suspended}
 \left\{\begin{array}{l}
 \dot X=F(X)+(\gamma(1-x)\sin(s),0,0)\\
 \dot s=2\omega\pmod{2\pi} .
\end{array} \right.
 \end{equation}
If $s_0\in \EU^1$, then the solutions of \eqref{eq:suspended} with initial condition $(0,0,\pm1,s_0)$ are $\pi/\omega$--periodic and we abuse notation denoting them by $\vv$ and $\ww$.

Under the condition that  $(\alpha-\beta)^2<4\alpha$ and $\beta \neq \alpha -2$,  the  first return map to a suitably rescaled section transverse  to the connection $[\vv\to\ww]$ was derived in \cite[Theorem 1]{LR18}.
 It is shown there
that it may be reduced to a two-dimensional 
map $G(y,s)$ defined on the cylinder
$$
\{(y,s):\  s\in\RR\pmod{2\pi},\quad y>0 \}
$$
and given by 
\begin{equation}
\label{map_G}
G (y,s)=
\left(y^{\delta^2}  +\gamma y^{\delta^2-\delta}\left(1+k\sin s\right),s-\frac{ \ln y}{K}\right),
\end{equation}
where  $K=(\alpha+\beta)^2/2\alpha>0$ and $k>0$  is a  constant depending on $\alpha$, $\beta$
and $\omega$. 
In the reduction of \cite{LR18} the variable $y$ is associated to $X$ while $s$ retains its meaning of perturbation time.
The explicit value of $k$ is difficult to estimate
 and this is why a full  description of the dynamics of \eqref{map_G} for all $k > 0$ is  welcome. The case $k = 0$ is also studied for the sake of completeness.

 The dynamics of maps similar to \eqref{map_G} has been studied in the literature in several contexts, either physical or biological, emphasising the role of different parameters.  
 In the context of the May--Leonard system,  the authors of \cite{DT3, TD1} have explored the cases $\delta   \gtrsim 1$ and $\delta\gg1$.  
 
 In the weakly attracting case ($\delta \gtrsim 1$), for the May--Leonard system
Dawes and Tsai found in \cite{DT3, TD1} saddle-node bifurcations and frequency locked solutions. 
Their proof is based on an equivalence to the model of a damped pendulum with constant torque. 
 If $\omega \ll 1$ an invariant torus has been found using the Afraimovich Annulus Principle.
 These results have been refined in \cite{LR18} where it is shown  that there exists
  a line of Bogdanov--Takens bifurcations in the space $(T, k, \gamma)$, where $k<1$ and  $T>0$ is the period of a periodic point of \eqref{map_G}.

     In the strongly attracting case ($\delta \gg 1$), the authors of  \cite{ DT3, TD1}   proved an equivalence to a family of circle maps if $\omega \gg 1$ and chaos when $\omega \approx 0$. These results have been refined in \cite{Rodrigues_SIAM} where it is shown  that,  for  values of the amplitude of the periodic forcing in a set of positive Lebesgue measure,  the dynamics is dominated by strange attractors with fully stochastic properties. The proof is performed by using the Wang and Young's theory of rank-one strange attractors \cite{WY2003}.   The attracting case with $k<1$ and $\omega \rightarrow +\infty$ was also discussed in \cite{LR23.1, TD1}. 
     
 Here we formalize what is meant by ``\emph{weakly}'' and ``\emph{strongly''} attracting cases and we extend the bifurcation analysis started in \cite{LR18} to the cases $k \geq 1$ and $k=0$. 
 Our focus is on  the analysis of the transitions between the different bifurcation diagrams according to the  parameters $\gamma$ and $\delta$.

 \subsection{The family of circle maps}\label{subsec:circleMap}

In dealing with the bifurcation and stability of periodic solutions of \eqref{general} we first introduce an auxiliary parameter $\tau$ and look for  fixed points of the map
$$
\GG(y,s)=G(y,s)+\left(0,\dfrac{\ln\tau}{K}\right),
$$ 
i.e., solutions of the equation
 \begin{equation}\label{eq:PeriodicSolutions}
\GG(y,s)=\left(
y^{\delta^2} +\gamma y^{\delta^2-\delta }\left(1+k\sin s\right),
s-\frac{\ln y}{K}\right)
+\left(0,\dfrac{\ln\tau}{K}\right)=(y,s)
\quad\mbox{for}\quad 0<\tau\le 1 .
\end{equation} 
Thus, the value $\tau=1$ corresponds to constant solutions of \eqref{general} and the behaviour in the limit $\tau\to 0$ describes solutions whose period tends to $\infty$.
 The auxiliary parameter $\tau$ will be removed  later to yield solutions of the original problem.
 
 The geometry of the set of fixed points of the map $\GG(y,s)$ is obtained here for different values of $\gamma>0$ and $\delta>1$.
 The solution of the second component of equation \eqref{eq:PeriodicSolutions} is  $y=\tau$,
 hence we look at zeros of the family of circle maps $s\mapsto g(\tau,s)$  where
 \begin{equation}\label{eq:g0}
 g(\tau,s)=\tau^{\delta^2} +\gamma \tau^{\delta^2-\delta }\left(1+k\sin s\right)-\tau .
 \end{equation}
The solutions of the equation $g(\tau,s)=0$ are drawn on a cylinder
$$
\CC=\{(\tau,s):\ 0<\tau\le 1,\quad s\in\EU^1\}.
$$
Note that although  the equation \eqref{eq:PeriodicSolutions} is not defined for $\tau=0$, the expression \eqref{eq:g0} satisfies $g(0,s)\equiv0$.

\section{Bifurcation diagrams for the circle maps }\label{sec:bifCylinder}

The equation $g(\tau,s)=0$ defines a curve on the cylinder $\CC$, the bifurcation diagram, that varies with the parameters $\gamma$, $\delta$ and $k$.

 \subsection{Preliminary results and  definitions}
 
 We start by establishing some standard terminology following the conventions in Golubitsky and Schaeffer \cite{vol1} and Keyfitz  \cite{Key}.
  \begin{definition}\label{def:bifDiagram}
 The {\em bifurcation diagram} for the map $g:\CC\seta\RR$ with {\em distinguished parameter}  $\tau$ is 
the  set  of solutions of $g(\tau,s)=0$ on  $\CC$, under the equivalence relation defined below.
 \end{definition}

 \begin{definition}\label{def:contactEq}
 The bifurcation diagrams for  the functions $g,h: \CC\rightarrow\RR$
 are {\em contact equivalent} if there is a diffeomorphism of the form $(T(\tau),S(\tau,s))$ and there is a function $A(\tau,s)>0$ such that 
 $$
 h(\tau,s)=A(\tau,s)g(T(\tau),S(\tau,s)),
 $$ 
 i.e., there is a  family of changes of coordinates  $S$  in the circle that preserves the fibration by the parameter $\tau$.
 \end{definition}
  \begin{definition}
The bifurcation diagram for a function $g$ is  {\em persistent} if any small  perturbation of  $g$ in the $C^1$--Whitney topology yields a contact equivalent bifurcation diagram. 
 \end{definition}
Bifurcation diagrams under Definition~\ref{def:contactEq} of contact equivalence have been classified by Keyfitz \cite{Key}, see also Golubitsky and Schaeffer \cite{vol1}.

 \begin{definition}
We say a point  $(\tau_0,s_0)\in \CC$ such that $g(\tau_0,s_0)=0$
is a {\em bifurcation point} when,  in any neighbourhood of $(\tau_0,s_0)$ in $\CC$,  the number $n(\tau)$ of solutions $s$ (for each $\tau$) to the equation $g(\tau,s)=0$
 changes at $\tau_0$.
  \end{definition}

  The only generic bifurcation points found in our  context are \emph{fold points} at some $(\tau,s)=(\tau_0,s_0)$,
 where $\dfrac{\partial g}{\partial s}(\tau_0,s_0)=0$.
We distinguish between \emph{supercritical folds} where locally there are two solutions to $g(\tau,s)=0$ for $\tau>\tau_0$, a single solution at $\tau_0$ and no solution for $\tau<\tau_0$;
\emph{subcritical folds} are those where the preceding inequalities are reversed.
More details can be found in  \cite{vol1,Key}.
 \bigbreak

Both expressions  \eqref{eq:PeriodicSolutions} and \eqref{eq:g0} depend on the parameters $\delta$ and $\gamma$.
Accordingly, we draw {\em transition diagrams}\footnote{These are usually also called {\em bifurcation diagrams} in the literature, we use a different terminology to avoid  confusion.}  where the $(\delta,\gamma)$ plane is divided by curves into regions corresponding to  contact equivalent persistent bifurcation diagrams. 
Typical  transitions occur either when $g(\tau,s)$ has a Morse critical point, or when the solution curve $g(\tau,s)=0$ runs into the boundary of $\CC$.
 We start with some more notation and a simple result that will be used in the sequel.

The equation \eqref{eq:g0} for $g(\tau,s)=0$ can be rewritten as
\begin{equation}\label{eq:Fdelta}
\gamma(1+k\sin s)=\tau^{-\delta^2+\delta+1}-\tau^\delta.
\end{equation}

Bifurcation diagrams are obtained using the properties of the right hand side of equation \eqref{eq:Fdelta}, described in the next lemma that is stated for reference without its elementary proof.
The \emph{golden number} $\Phi=(1+\sqrt{5})/2$, a solution of $p(\delta)=0$ of
\begin{equation}\label{eq:pdelta}
p(\delta)=-\delta^2+\delta+1,
\end{equation}
plays a crucial role in this analysis.

\begin{figure}[hht]
\parbox{50mm}{\begin{center}
\includegraphics[width=45mm]{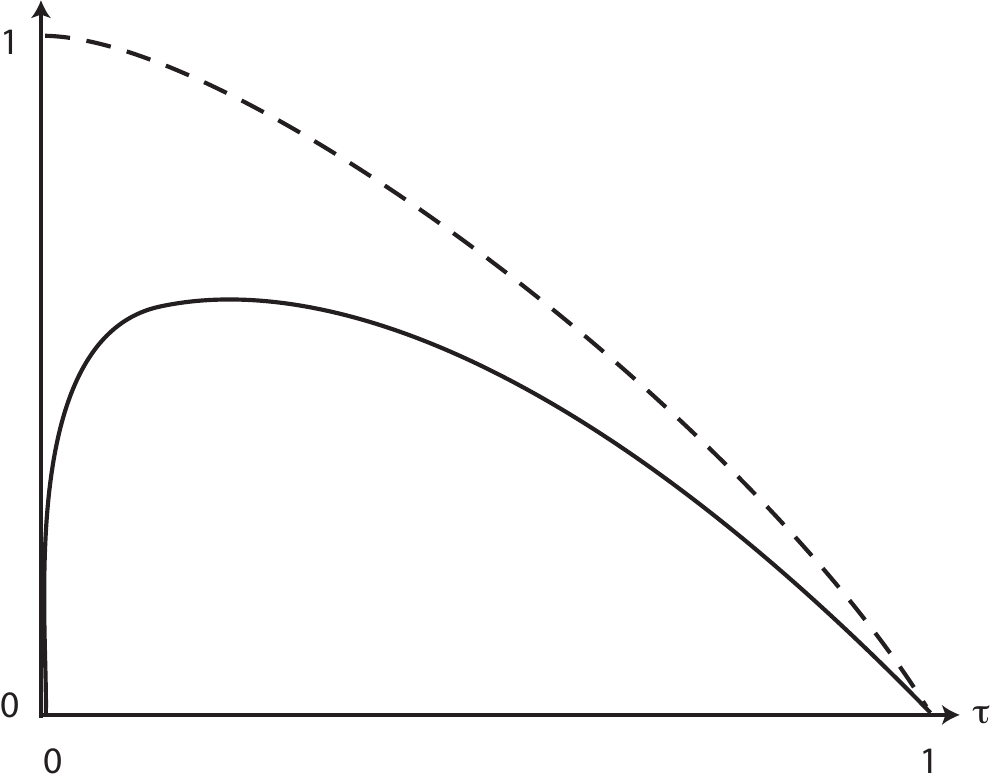}\\
$1<\delta\le\Phi$
\end{center}}  \qquad
\parbox{50mm}{\begin{center}
\includegraphics[width=45mm]{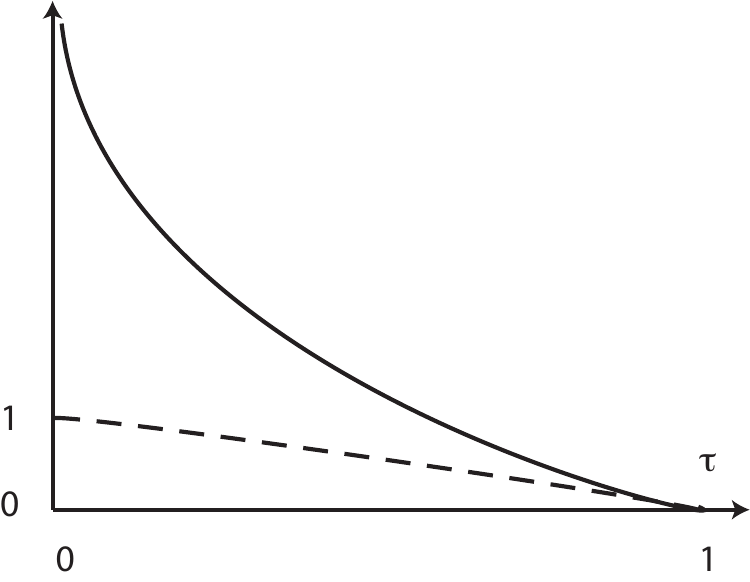}\\
 $\delta\ge \Phi$ 
\end{center}}\\
\caption{\small Graph of $\tau\to\FF_\delta(\tau)$ for different values of $\delta$.
 On the left the solid line is the graph for $1<\delta<\Phi$, the dashed line  is the graph of $\FF_\Phi(\tau)$.
On the right  the
solid line is the graph for
$\delta>\Phi$, the dashed line is the graph of $\FF_\Phi(\tau)$.}
\label{fig:Fdelta}
\end{figure}

\begin{lemma}\label{lem:Fdelta}
The following properties, illustrated in Figure~\ref{fig:Fdelta}, hold for the map $\FF_\delta:(0,1]\seta\RR$  given by
$\FF_\delta(\tau)=\tau^{p(\delta)}-\tau^\delta$:
\begin{enumerate}
\item\label{item:F>0}
$\FF_\delta(\tau)>0$ for all $\tau\in(0,1)$  and $\FF_\delta(1)=0$;
\item\label{item:Fdelta=phi}
$\FF_\Phi(\tau)$ is monotonically decreasing and $\dpt\lim_{\tau\to 0}\FF_\Phi(\tau)=1$;\item\label{item:Fdelta>phi}
if $\delta>\Phi$, then $\FF_\delta(\tau)$ is monotonically decreasing, with $\dpt\lim_{\tau\to 0^+}\FF_\delta(\tau)=\infty$;
\item\label{item:Fdelta<phi}
if $1<\delta<\Phi$, then  $\dpt\lim_{\tau\to 0}\FF_\delta(\tau)=0$ and $\FF_\delta(\tau)$ has a global maximum $M_\FF(\delta)>0$ at 
$$
\tau_m=\left(\dfrac{p(\delta)}{\delta} \right)^{1/(\delta^2-1)}
\mbox{with}\qquad
 M_\FF(\delta):= \left(\dfrac{p(\delta)}{\delta} \right)^{\frac{p(\delta)}{\delta^2-1}}- 
 \left(\dfrac{p(\delta)}{\delta} \right)^{\frac{\delta}{\delta^2-1}};
$$
\item  $\dpt \lim_{\delta \rightarrow \Phi^-}\tau_m=0$,
 $\dpt \lim_{\delta \rightarrow \Phi^-}M_\FF(\delta) =1$ and $\dpt \lim_{\delta \rightarrow 0 }M_\FF(\delta) =0$.
\end{enumerate}
\end{lemma}

\begin{figure}[hht]
 \includegraphics[width=135mm]{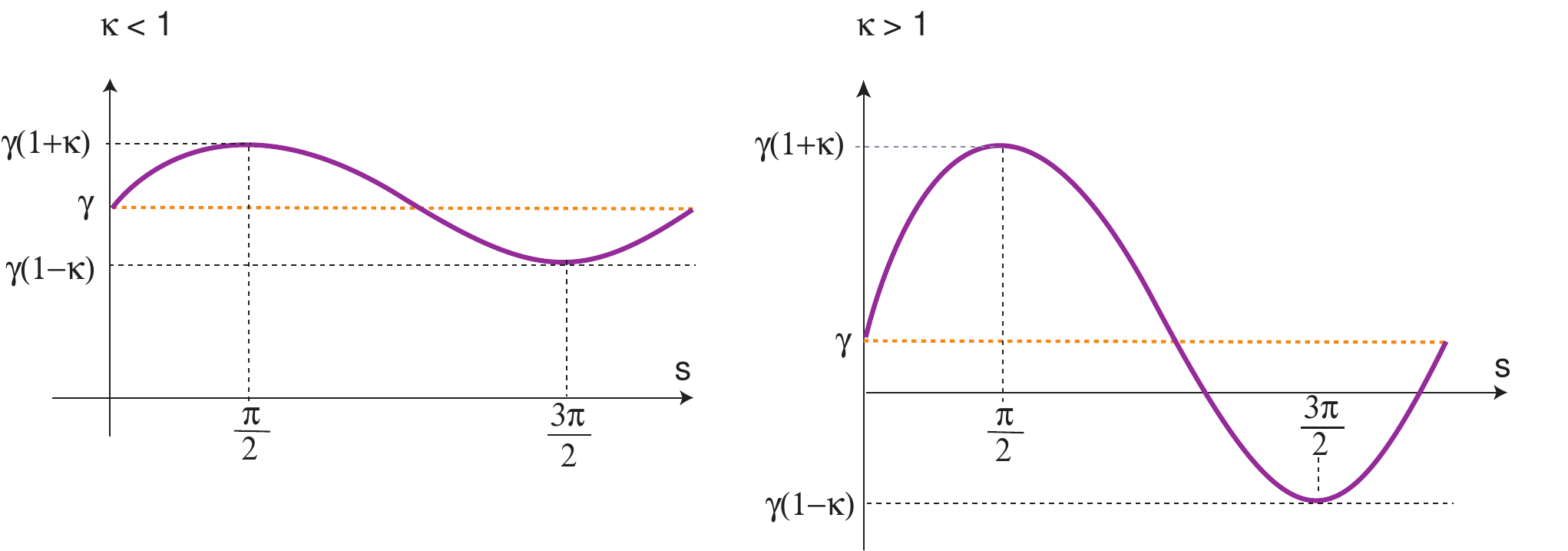} 
\caption{\small Graph of $s\mapsto \gamma(1+k\sin s)$ for  $k<1$ (left) and $k>1$ (right).}
\label{fig:graph_phi1_new}
\end{figure}

 Since $\FF_\delta(\tau)>0$, only the positive values of the map  $s\mapsto \gamma(1+k\sin s)$ are of interest for the bifurcation diagram of $g$.
Thus for $\gamma>0$  the different signs of $k-1$ give rise to  different diagrams as can be deduced from the graphs in   Figure \ref{fig:graph_phi1_new}. 
 We consider separately the cases $0<k<1$ and $k>1$ starting with the simpler case $k>1$.
After these, we describe  the limit case $k=0$ that corresponds to an autonomous perturbation of equations \eqref{general} and the transition case $k=1$.

\subsection{Case $k>1$}\label{subsec:k>1}

The following result establishes the  $(\tau,s)$ bifurcation diagrams  associated to  the map \eqref{eq:g0} as well as the transitions in the $(\delta, \gamma)$ plane for $k>1$. 
 We recommend that the reader follows the proof in a cinematic way: 
 fix a  value of  $\tau$,
 find the value of $\FF_\delta(\tau)$ and then intersect the horizontal line $y=\FF_\delta(\tau)$ (which is moving up and/or down) with the graph of $s\mapsto\gamma(1+k\sin s)$.
The bifurcation diagrams always consist of regular curves in $\CC$. 
For $k>1$ these are open curves, i.e., curves homeomorphic to an open interval in $\RR$.

\begin{figure}[hht]
\parbox{50mm}{\begin{center}
\includegraphics[width=35mm]{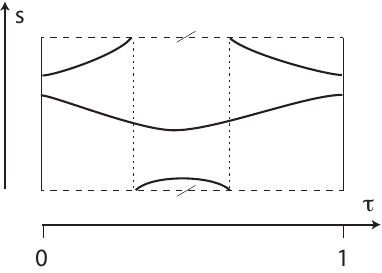}\\
bifurcation diagram\\
region A 
\end{center}}\quad
\parbox{50mm}{\begin{center}
\includegraphics[width=35mm]{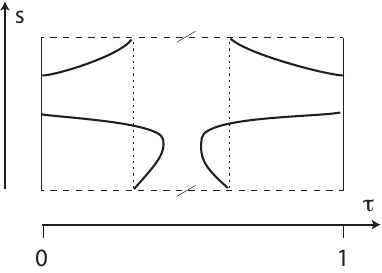}\\
bifurcation diagram\\
region B
\end{center}}\quad
\parbox{50mm}{\begin{center}
\includegraphics[width=35mm]{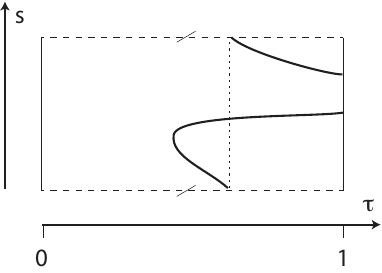}\\
bifurcation diagram\\
region C
\end{center}}\\
\parbox{50mm}{\begin{center}
\includegraphics[width=39mm]{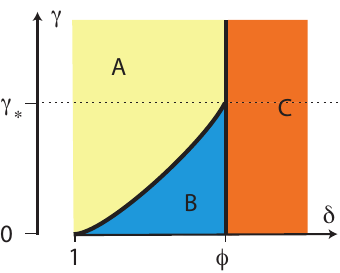}\\
transition diagram\\
$\gamma_*=1/(1+k)$
\end{center}}\quad
\parbox{50mm}{\begin{center}
\includegraphics[width=35mm]{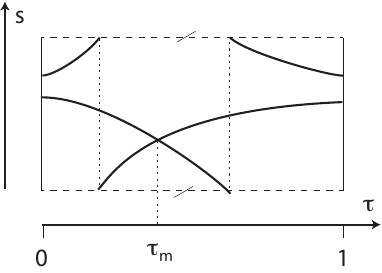}\\
bifurcation diagram\\
transition A $\leftrightarrow$ B\\
$\gamma=M_\FF(\delta)/(1+k)=:\gamma_+(\delta)$
\end{center}}\quad
\parbox{50mm}{\begin{center}
\includegraphics[width=35mm]{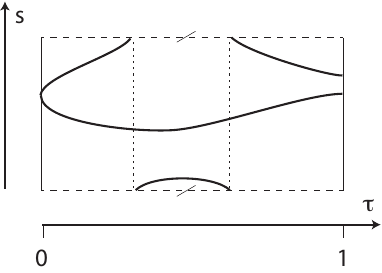}\\
bifurcation diagram\\
transition A $\leftrightarrow$ C
\end{center}}
\caption{\small Illustration of Theorem \ref{th:k>1}. Top: persistent bifurcation diagrams on the cylinder with $k>1$, labelling of regions as in Theorem~\ref{th:k>1}.
Bottom:  
transition diagram for the map \eqref{eq:g0} with $k>1$ and bifurcation diagrams  for parameters on the curves separating the regions.
The horizontal boundaries in the bifurcation diagrams are identified, we indicate this by the line /.}
\label{fig:k>1}
\end{figure}

\begin{theorem}\label{th:k>1}
For $k>1$ the persistent bifurcation diagrams of the map $g$, shown in Figure~\ref{fig:k>1}, satisfy:
\begin{enumerate}
\renewcommand{\theenumi}{(\Alph{enumi})}
\renewcommand{\labelenumi}{{\theenumi}}
\item\label{item:A}
if $\delta<\Phi$ and $\gamma> M_\FF(\delta)/(1+k)=:\gamma_+(\delta)$, then  the diagram comprises two open curves along  the length of the cylinder with no fold points;
\item\label{item:B}
if $\delta<\Phi$ and $0<\gamma< M_\FF(\delta)/(1+k):=\gamma_+(\delta)$, then the diagram comprises two open curves not covering the length  of the cylinder, one with a subcritical, the other with a supercritical fold point;
\item\label{item:C}
if $\delta>\Phi$, then the diagram is an open curve not covering the length  of the cylinder with a supercritical fold point;
\setcounter{lixo}{\value{enumi}}
\end{enumerate}
where $M_\FF(\delta)>0$ is defined  in assertion \eqref{item:Fdelta<phi} of Lemma~\ref{lem:Fdelta}.
\end{theorem}

\begin{proof}
The left hand side $s\mapsto \gamma(1+k\sin s)$ of \eqref{eq:Fdelta} lies between 
$\gamma(1-k)<0<\gamma(1+k)$, where $\gamma(1-k)<0$.
Since $\FF_\delta(\tau)>0$ for $\tau\in(0,1)$ and  $k>1$,   then solutions of equation \eqref{eq:Fdelta} only exist when  $\sin s>-1/k$.
Solutions of equation \eqref{eq:Fdelta}  lie in the intersection of the graph of $s\mapsto \gamma(1+k\sin s)$ with the horizontal line at height   $\FF_\delta(\tau)$. 
Therefore the bifurcation diagrams can be obtained observing the position of the graph of $z=\FF_\delta(\tau)$ with respect to the horizontal line $z=\gamma(1+ k)$.
Fold points occur when $\FF_\delta(\tau)=\gamma(1+k)$, they  are subcritical when  $\FF'_\delta(\tau)>0$, supercritical when $\FF'_\delta(\tau)<0$.

If $\delta<\Phi$ and $\gamma>M_\FF(\delta)/(1+k)$, then $\FF_\delta(\tau)<\gamma(1+k)$ for all $\tau\in(0,1]$.
For each $\tau$ there are two values of $s\in [0, 2\pi)$ satisfying equation  \eqref{eq:Fdelta} as in \ref{item:A}.

If $\delta<\Phi$ and $M_\FF(\delta)/(1+k)>\gamma$, then $\FF_\delta(\tau)=\gamma(1+k)$ at two values of $\tau\in(0,1]$ that correspond to a subcritical and a supercritical fold.
For $\tau$ between these values there are no solutions to \eqref{eq:Fdelta}.
If $\tau$ does not lie in  this interval there are  two values of $s$ satisfying equation  \eqref{eq:Fdelta} as in \ref{item:B}. 

For $\delta>\Phi$ the map $\FF_\delta(\tau)$ is monotonically decreasing and hence there is a supercritical fold at the solution $\tau_0\in (0, 1]$ of $\FF_\delta(\tau_0)=\gamma(1+k)$ meaning that for each $\tau>\tau_0$ there are two  $s\in [0, 2\pi)$ yielding solutions $(s,\tau)$
of equation \eqref{eq:Fdelta} as in \ref{item:C}. 
\end{proof}

\begin{figure}[hht]
\parbox{50mm}{
\begin{center}
\includegraphics[width=49mm]{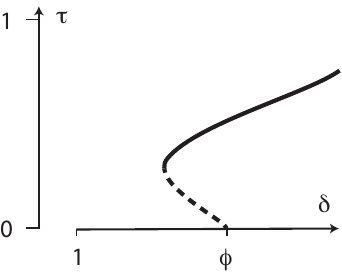}\\
$\gamma<1/(1+k)$
\end{center}} \quad
\parbox{50mm}{\begin{center}
\includegraphics[width=49mm]{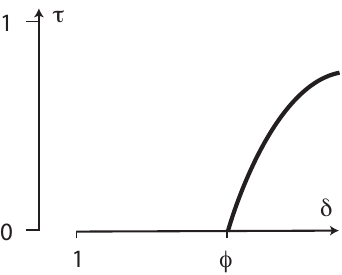}\\
$\gamma>1/(1+k)$
\end{center}}\\
\caption{\small Fold points for the map \eqref{eq:g0} with $k>1$. Solid lines are supercritical folds, dashed lines are subcritical. }
\label{fig:Foldsk>1}
\end{figure}

If $k>1$
the fold points correspond to $s=\pi/2$; their position as $\delta$ varies is shown in Figure~\ref{fig:Foldsk>1}.

The bifurcation diagrams on the boundary lines of regions \ref{item:A}--\ref{item:C} of Theorem~\ref{th:k>1} are  shown in Figure~\ref{fig:k>1} at the bottom. 
They represent transitions between the regions of persistent behaviour as follows:
\begin{itemize}
\item\label{item:AB} 
If $\delta<\Phi$ and  $\gamma=M_\FF(\delta)/(1+k) =\gamma_+(\delta)$, then the two  fold points on the curves that exist for smaller $\gamma$ coalesce. This transition is locally described by the level curves of a Morse function around a saddle point.

\item\label{item:transitionBPhi}
In the limit $\delta\to\Phi^-$ with $\gamma<\gamma_*=1/(1+k)$ 
 the whole branch that contains a subcritical 
fold collapses into the boundary $\tau=0$ of $\CC$.
This is the transition \ref{item:B}--\ref{item:C}.

\item\label{item:transitionAPhi}
In the limit $\delta\to\Phi^-$ with $\gamma>\gamma_*=1/(1+k)$ 
 the two branches of the bifurcation diagram come together at a supercritical  fold point on the boundary $\tau=0$ of the cylinder $\CC$ --- see Figure~\ref{fig:Foldsk>1}.
As $\delta$ increases from $\Phi$ the fold point moves into larger values of $\tau$.
This is the non standard transition  \ref{item:A}--\ref{item:C}.
\end{itemize}

\subsection{Case $0<k<1$}

\begin{figure}[hht]
\parbox{50mm}{\begin{center}
\includegraphics[width=49mm]{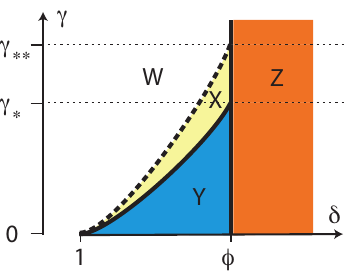}\\
transition diagram\\
$\gamma_*=1/(1+k)$\\
$\gamma_{**}=1/(1-k)$
\end{center}}\\
\parbox{50mm}{\begin{center}
\includegraphics[width=35mm]{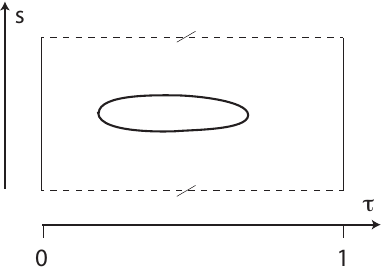}\\
bifurcation diagram\\
region X
\end{center}}\quad
\parbox{50mm}{\begin{center}
\includegraphics[width=35mm]{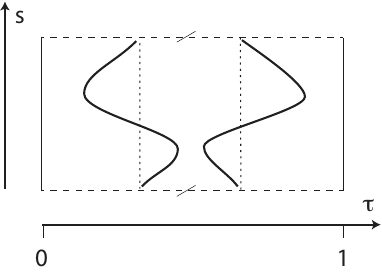}\\
bifurcation diagram\\
region Y
\end{center}}\quad
\parbox{50mm}{\begin{center}
\includegraphics[width=35mm]{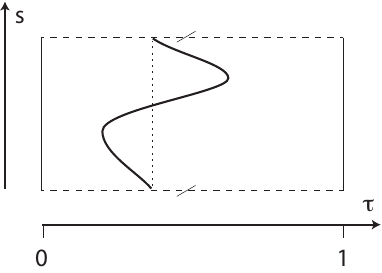}\\
bifurcation diagram\\
region Z
\end{center}}\\
\parbox{50mm}{\begin{center}
\includegraphics[width=35mm]{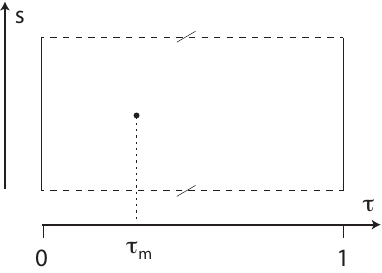}\\
bifurcation diagram\\
transition W $\leftrightarrow$ X\\
$\gamma=M_\FF(\delta)/(1-k) =:\gamma_-(\delta)$
\end{center}}\quad
\parbox{50mm}{\begin{center}
\includegraphics[width=35mm]{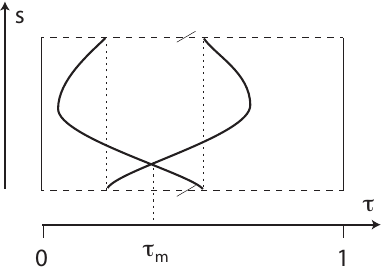}\\
bifurcation diagram\\
transition X $\leftrightarrow$ Y\\
$\gamma=M_\FF(\delta)/(1+k) =:\gamma_+(\delta)$
\end{center}}\quad
\parbox{50mm}{\begin{center}
\includegraphics[width=35mm]{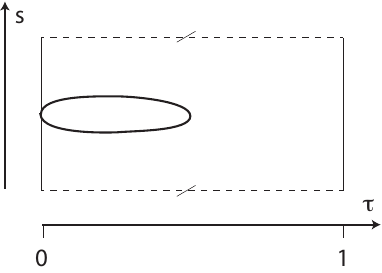}\\
bifurcation diagram\\
transition X $\leftrightarrow$ Z
\end{center}}\\
\caption{\small Illustration of Theorem~\ref{th:k<1}. Top: transition diagram for the map $g$ defined in \eqref{eq:g0} with $0<k<1$. 
In the region W there are no  zeros of the map $g$.
Center: persistent bifurcation diagrams.
Bottom:  bifurcation diagrams for parameters in the curves separating the regions.
For $\delta=\Phi$ a branch of the diagram is absorbed by the boundary $\tau=0$.
The horizontal boundaries in the bifurcation diagrams are identified, we indicate this by the line /.}
\label{fig:k<1}
\end{figure}

The next result deals with the 
$(\tau,s)$ bifurcation diagrams  associated to the map \eqref{eq:g0} 
 for $k<1$ as well as the conditions on $\delta$ and $\gamma$ under which they occur.
 
\begin{theorem}\label{th:k<1}
For $0<k<1$ the persistent bifurcation diagrams of $g$,
shown in Figure~\ref{fig:k<1}, satisfy:
\begin{enumerate}
\renewcommand{\theenumi}{(\Alph{enumi})}
\renewcommand{\labelenumi}{{\theenumi}}
\setcounter{enumi}{22}
\item\label{item:W}
if $\delta<\Phi$ and $\gamma> M_\FF(\delta)/(1-k) =:\gamma_-(\delta)$, then  $g(\tau,s)\ne 0$ at all points
in the cylinder $\CC$;
\item\label{item:X}
if $\delta<\Phi$ and $ \gamma_+(\delta)=:M_\FF(\delta)/(1+k)<\gamma< M_\FF(\delta)/(1-k)=:\gamma_-(\delta)$, then the diagram is a contractible closed curve in the cylinder $\CC$ with two fold points;
\item\label{item:Y}
if $\delta<\Phi$ and $0<\gamma< M_\FF(\delta)/(1+k)=:\gamma_+(\delta)$, then the diagram comprises two noncontractible closed curves in the cylinder $\CC$ with two fold points each;
\item\label{item:Z}
if $\delta>\Phi$, then the diagram is a noncontractible closed curve in the cylinder $\CC$ with two fold points;
\setcounter{lixo}{\value{enumi}}
\end{enumerate}
where $M_\FF(\delta)>0$ is defined  in assertion \eqref{item:Fdelta<phi} of Lemma~\ref{lem:Fdelta}.
\end{theorem}

\begin{proof}
Since $0<k<1$, then  the left hand side of equation \eqref{eq:Fdelta} varies in the range $[\gamma(1-k), \gamma(1+k)]$, where $0<\gamma(1-k)<\gamma(1+k)$.
Fold points occur when the graph of $\FF_\delta(\tau)$ attains these limiting values.
Fold points are subcritical when either $\FF_\delta(\tau)=\gamma(1+k)$ and $\FF'_\delta(\tau)>0$ or when $\FF_\delta(\tau)=\gamma(1-k)$ and $\FF'_\delta(\tau)<0$.
Reversing the signs of the derivatives yields the conditions for supercritical folds.
Hence, the bifurcation diagrams can be obtained observing the position of the graph of $z=\FF_\delta(\tau)$ with respect to the horizontal lines $z=\gamma(1\pm k)$.
This establishes the result as follows.

If $\delta<\Phi$ and $M_\FF(\delta)/(1-k)<\gamma$, then $\FF_\delta(\tau)$ never takes values in the interval $[\gamma(1-k),\gamma(1+k)]$.
The bifurcation diagram is empty establishing \ref{item:W}.

If $\delta<\Phi$ and  $M_\FF(\delta)/(1+k)<\gamma<M_\FF(\delta)/(1-k)$, then as $\tau$ increases there is first a supercritical fold followed by a subcritical one. 
Solutions of equation \eqref{eq:Fdelta} only exist for $\tau$ in the interval between the two fold points.
The diagram is  a contractible curve on the cylinder that indeed reduces to a point at  $\gamma=M_\FF(\delta)/(1+k)$ proving assertion \ref{item:Z}.

Still in the case $\delta<\Phi$, if $\gamma<M_\FF(\delta)/(1+k)$, then the graph of $\FF_\delta(\tau)$ crosses twice each one of the two limiting values $\gamma(1\pm k)$ yielding four fold points.  
For $\tau$ in the interval between the two points where $\FF_\delta(\tau)=\gamma(1-k)$  there are no solutions of equation \eqref{eq:Fdelta}.
The bifurcation diagram consists of two disjoint  curves that go round the cylinder each one containing two fold points, as in \ref{item:Y}.

Finally, when $\delta>\Phi$, then $\FF'_\delta(\tau)<0$ and the graph of $\FF_\delta(\tau)$ crosses each  horizontal line $z=\gamma(1\pm k)$ once.
The bifurcation diagram consists of a  curve that goes round the cylinder containing two fold points, as in \ref{item:Z}.
\end{proof}

\begin{figure}[hht]
\parbox{50mm}{\begin{center}
\includegraphics[width=49mm]{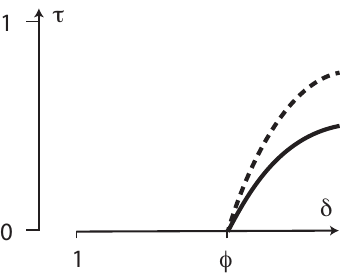}\\
$\gamma>1/(1-k)$
\end{center}}\quad
\parbox{50mm}{\begin{center}
\includegraphics[width=49mm]{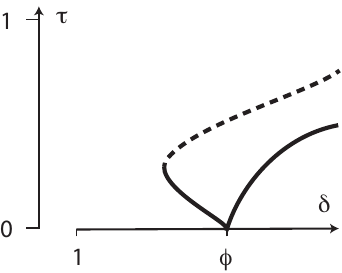}\\
$1/(1+k)<\gamma<1/(1-k)$
\end{center}}\quad
\parbox{50mm}{\begin{center}
\includegraphics[width=49mm]{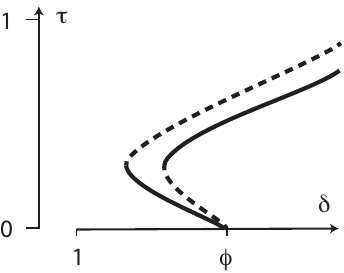}\\
$\gamma<1/(1+k)$
\end{center}}\\
\caption{\small Fold points for \eqref{eq:g0} with $k<1$. Solid lines are supercritical folds, dashed lines are subcritical.}
\label{fig:Foldsk<1}
\end{figure}

The position of fold points in the $(\delta,\tau)$ plane is shown in Figure~\ref{fig:Foldsk<1} for $k<1$.
In the transition between the regions described in \ref{item:W}--\ref{item:Z} 
of Theorem~\ref{th:k<1}
 the behaviour is as shown in the bottom part of Figure~\ref{fig:k<1}.
 There are  three  types of transition: generic codimension 1 transitions between regions \ref{item:W}, \ref{item:X} and \ref{item:Y},
 transitions where a fold point moves to the boundary of $\CC$,
  and  global transitions at 
 the boundary of the cylinder at $\delta=\Phi$ where a whole curve collapses into the boundary of $\CC$.
 More specifically, they are:
\begin{itemize}
\item\label{item:XY} 
if $\delta<\Phi$ and  $\gamma=M_\FF(\delta)/(1+k)=\gamma_+(\delta)$, then two of the fold points of the curves that exist for smaller $\gamma$ coalesce. This transition is locally described by the level curves of a Morse function around a saddle point, this is the transition between cases \ref{item:X} and \ref{item:Y};
\item\label{item:XW} 
if $\delta<\Phi$ and  $\gamma=M_\FF(\delta)/(1-k) =\gamma_-(\delta)$, then the contractible  curve that exists for smaller $\gamma$ shrinks to a point. This transition is locally described by the level curves of a Morse function around a minimum in the transition  \ref{item:X}--\ref{item:W};
\item\label{item:transitionPhiWZ}
in the limit $\delta\to\Phi^+$ with $\gamma>M_\FF(\delta)/(1-k) =\gamma_-(\delta)$ the noncontractible closed curve in the bifurcation diagram collapses into the boundary $\tau=0$ of $\CC$.
This is the non standard transition \ref{item:W}--\ref{item:Z};
\item\label{item:transitionPhiXZ}
in the limit $\delta\to\Phi^-$ with $ \gamma_+(\delta)=M_\FF(\delta)/(1+k)<\gamma<M_\FF(\delta)/(1-k) =\gamma_-(\delta)$ the closed curve in the bifurcation diagram touches the boundary $\tau=0$ of $\CC$ turning into a noncontractible curve on the cylinder.
This is the non standard transition \ref{item:X}--\ref{item:Z},
 a fold point at the boundary of $\CC$;
\item\label{item:transitionPhiYZ}
in the limit $\delta\to\Phi^-$ with $\gamma<M_\FF(\delta)/(1+k) =\gamma_+(\delta)$ a branch of the bifurcation diagram collapses into the boundary $\tau=0$ of $\CC$.
This is the non standard transition \ref{item:Y}--\ref{item:Z}.
\end{itemize}

\subsection{Transition case $k=1$}
When $k<1$ increases to the value 1,  fold points in the bifurcation diagrams move to the components  $\tau=0$ and $\tau=1$ of the boundary of $\CC$. 
The region \ref{item:W} studied in Theorem \ref{th:k<1} disappears since $(1-k)\gamma=0<M_\FF(\delta)$ for all $\delta>1$.
The next result is a summary of this transition.
Taking into account that for $k=1$ one gets $ M_\FF(\delta)/(1+k)=M_\FF(\delta)/2$. 
The proof  is analogous to that of Theorems~\ref{th:k>1} and \ref{th:k<1} and is omitted.

\begin{figure}[hht]
\parbox{50mm}{\begin{center}
\includegraphics[width=49mm]{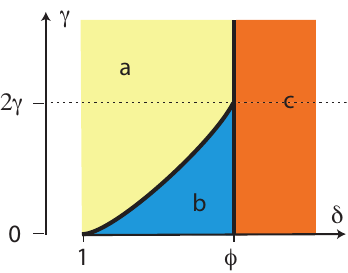}\\
transition diagram
\end{center}}\\
\parbox{50mm}{\begin{center}
\includegraphics[width=35mm]{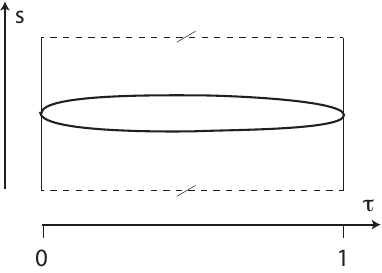}\\
bifurcation diagram\\
region a
\end{center}}\quad
\parbox{50mm}{\begin{center}
\includegraphics[width=35mm]{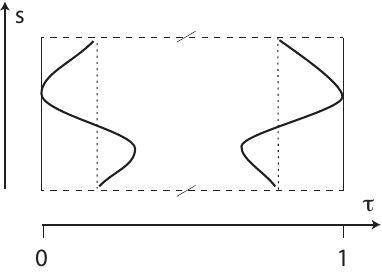}\\
bifurcation diagram\\
region b
\end{center}}\quad
\parbox{50mm}{\begin{center}
\includegraphics[width=35mm]{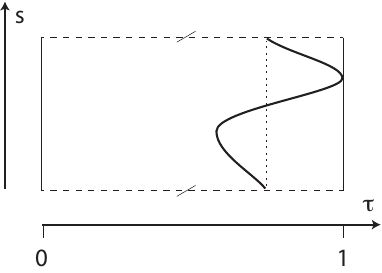}\\
bifurcation diagram\\
region c
\end{center}}
\caption{\small  Illustration of Theorem~\ref{th:k=1}. Top: transition diagram for  the map $g$ defined in \eqref{eq:g0} for the transition case $k=1$.  At the transition diagram on top, the curve separating regions \ref{item:a} and \ref{item:b} is the graph of $ \gamma(\delta)=M_\FF(\delta)/2$ for $\delta \in [1, \Phi)$.
Bottom: Persistent bifurcation diagrams for the map defined in  \eqref{eq:g0} with $k=1$. }
\label{fig:k=1}
\end{figure}

\begin{theorem}\label{th:k=1}
For $k=1$ the transition bifurcation diagrams of  $g$,
shown in Figure~\ref{fig:k=1}, satisfy:
\begin{enumerate}
\renewcommand{\theenumi}{(\alph{enumi})}
\renewcommand{\labelenumi}{{\theenumi}}
\item\label{item:a}
if $\delta<\Phi$ and $ M_\FF(\delta)/2<\gamma$,  then the diagram is a contractible closed curve in the cylinder with two fold points at $\tau=0$ and $\tau=1$;
\item\label{item:b}
if $\delta<\Phi$ and $0<\gamma< M_\FF(\delta)/2$,  then the diagram comprises two noncontractible closed curves in the cylinder with two fold points each where one of the curves has a fold point at $\tau=0$ and the other at $\tau=1$;
\item\label{item:c}
if $\delta>\Phi$, then the diagram is a noncontractible closed curve in the cylinder with two fold points one of which lies at $\tau=1$;
\setcounter{lixo}{\value{enumi}}
\end{enumerate}
 where $M_\FF(\delta)>0$ is defined  in  assertion \eqref{item:Fdelta<phi} of Lemma~\ref{lem:Fdelta}.
\end{theorem}

\subsection{Limit case $k=0$}
 Equation \eqref{eq:Fdelta} for $k=0$ becomes 
$\gamma =\tau^{-\delta^2+\delta+1}-\tau^\delta$ or equivalently $\gamma= \FF_\delta(\tau)$.
 Zeros of the map $g(\tau,s)$ correspond to the intersection of the graph of $\FF_\delta(\tau)$ with a horizontal line at height $\gamma$.

\begin{theorem}\label{th:k=0}
In the limit case $k=0$ the bifurcation diagrams of $g$
consist of  circles in $\CC$ with constant $\tau$ as follows:
\begin{enumerate}
\renewcommand{\theenumi}{(\roman{enumi})}
\renewcommand{\labelenumi}{{\theenumi}}
\item\label{item:k0i}
if $1<\delta<\Phi$ and $\gamma< M_\FF(\delta)$, then   $g(\tau,s)\ne 0$ at all points in $\CC$;
\item\label{item:k0ii}
if $1<\delta<\Phi$ and $\gamma> M_\FF(\delta) $, then the diagram consists of two circles in $\CC$;
\item\label{item:k0iii}
if $\delta>\Phi$,  then the diagram consists of one circle in $\CC$;
\setcounter{lixo}{\value{enumi}}
\end{enumerate}
 where $M_\FF(\delta)>0$ is defined  in  assertion \eqref{item:Fdelta<phi} of Lemma~\ref{lem:Fdelta}.
Moreover, in the transitions between these regions:
\begin{enumerate}
\renewcommand{\theenumi}{(\roman{enumi})}
\renewcommand{\labelenumi}{{\theenumi}}
\setcounter{enumi}{\value{lixo}}
\item\label{item:k0iv}
if $\delta=\Phi$ and $\gamma> M_\FF(\delta)=1 $, then  $g(\tau,s)\ne 0$ at all points in $\CC$;
\item\label{item:k0v}
if $\delta=\Phi$ and $\gamma\le M_\FF(\delta) =1 $, then the diagram consists of one circle in $\CC$;
\item\label{item:k0vi}
if $1<\delta<\Phi$ and $\gamma= M_\FF(\delta) $, then the diagram consists of one circle in $\CC$.
\end{enumerate}
\end{theorem}

\section{Stability and further bifurcations}\label{sec:stability}
In this section we study the  stability of the fixed points found in the previous section and we locate a codimension two bifurcation that has consequences for the dynamics.

The stability of  fixed points of the map $\GG$ (see \eqref{eq:PeriodicSolutions}) is obtained from the eigenvalues of its derivative
given by
\begin{equation}\label{eq:DG}
D\GG{} (y,s)= 
\left(
\begin{array}{ll} \delta^2   y^{\delta^2  -1} +\gamma(\delta^2-\delta)y^{-p(\delta)}\left( 1+k\sin s\right)&
 \gamma k y^{\delta^2-\delta} \cos s\quad\\
&\\
\dpt-\frac{1}{K y} &
1
 \end{array} \right) 
\end{equation}
where $p(\delta)$ was defined in \eqref{eq:pdelta} and $(y, s)\in \CC$.
They are easier to compute at the fold points, also called {\em saddle-nodes}.
\subsection{Stability near saddle-nodes}
\label{ss: SN}
In Section \ref{sec:bifCylinder} we have located generic fold points  through the intersection of graphs of smooth real valued functions defined in an interval. 
 In terms of dynamics the fold points correspond to saddle-nodes that, when they exist, occur
 at $s=\pi/2$ and $s=3\pi/2$.
In Proposition~\ref{prop:stabFold} below  we   identify the    stability of the  two emerging solution branches in some cases.
First we state a technical lemma that will be useful in the sequel whose proof is straightforward  and  left to the readers.

\begin{figure}[hht]
\includegraphics[width=80mm]{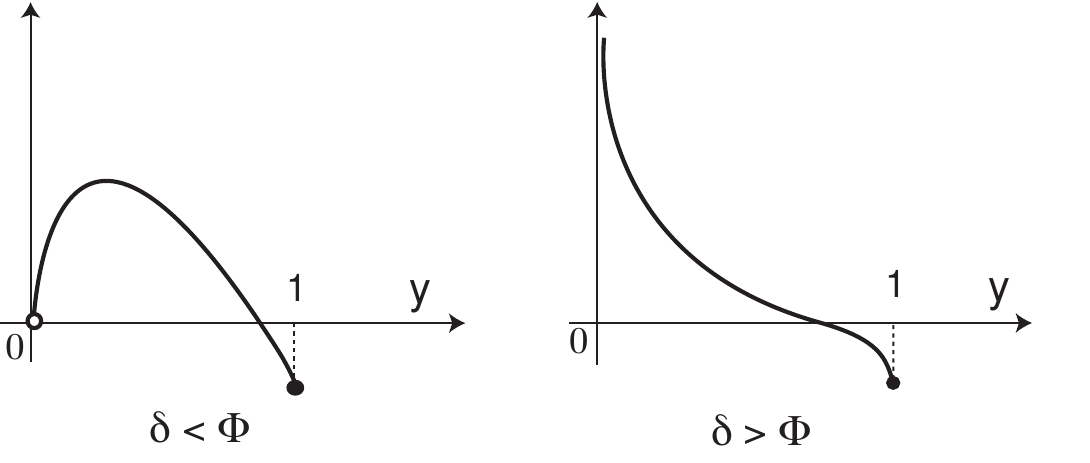}
\caption{\small Graphs of  the function $h_\delta(y)$ of Lemma~\ref{Lemma_aux} for different values of $\delta$.}
\label{fig:h}
\end{figure}

\begin{lemma}
\label{Lemma_aux}
 For $\delta>1$ and $p(\delta)$ defined in \eqref{eq:pdelta} the map $h_\delta:(0,1]\rightarrow \RR$ defined by 
 $$
 h_\delta(y)=  \dfrac{y^{p(\delta)}-\delta^2y^\delta}{\delta^2-\delta}
 $$
 is smooth and has a unique zero in its domain. Moreover the following assertions, illustrated in Figure~\ref{fig:h}, are true:
 \begin{enumerate}
 \item  \label{item:h1}
 $h_\delta(1)=-\dfrac{1+\delta}{\delta}<0$;
 \item \label{item:h2}
  $h_\delta(y)=0$ if and only if $y=\delta^{2/(1-\delta^2)} \in(0,1)$;
 \item  \label{item:h3}
 if  $\delta>\Phi$, then $p(\delta)<0$ and   $\dpt\lim_{y\to 0}h_\delta(y)=+\infty$;
 \item   \label{item:h4}
  if  $\delta<\Phi$, then $p(\delta)>0$ and   $\dpt\lim_{y\to 0}h_\delta(y)=0$.
  Moreover, 
 $\dfrac{dh_\delta}{dy}(y)=y^{\delta-1}\left(p(\delta)y^{1-\delta^2}-\delta^3\right)/(\delta^2-\delta)$
  with $\dpt\lim_{y\to 0}\dfrac{dh_\delta}{dy}(y)=+\infty$. 
  Therefore, if $0<y<\delta^{2/(1-\delta^2)}$, then $h_\delta(y)>0$.\\
 
 \end{enumerate}
\end{lemma}

\begin{figure}[hht]
\parbox{55mm}{\begin{center}
\includegraphics[width=45mm]{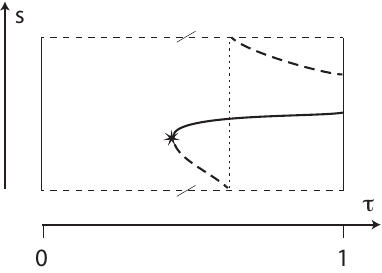}\\
region  \ref{item:C}\\
$k>1$\quad
$\delta>\Phi$ \\ $0<\gamma< M_\FF(\delta)/(1+k)=\gamma_+(\delta)$
\end{center}}\quad
\parbox{55mm}{\begin{center}
\includegraphics[width=45mm]{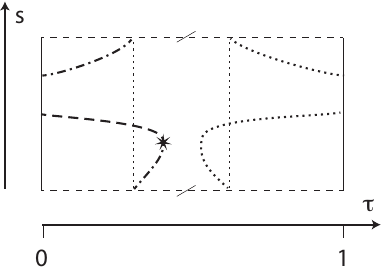}\\
region  \ref{item:B}\\
$k>1$\quad
$\delta<\Phi$
\end{center}}\\
\smallbreak
\parbox{55mm}{\begin{center}
\includegraphics[width=45mm]{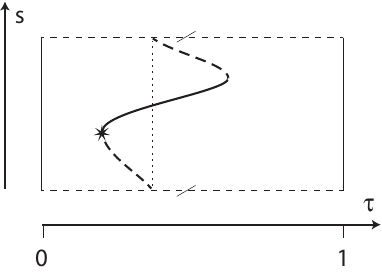}\\
region \ref{item:Z}\\
$0<k<1$\quad
$\delta>\Phi$
\end{center}}
\quad
\parbox{55mm}{\begin{center}
\includegraphics[width=45mm]{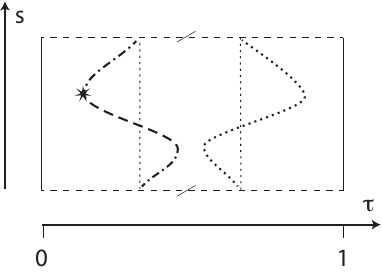}\\
region  \ref{item:Y}\\
$0<k<1$\quad
$\delta<\Phi$ \\ $0<\gamma< M_\FF(\delta)/(1+k) =\gamma_+(\delta)$
\end{center}}
\caption{\small Stability near fold points in the bifurcation diagrams of Proposition~\ref{prop:stabFold} assuming stability does not change away from the neighbourhood of the fold where it was computed, indicated as $\star$. Conventions: solid lines are attractors, dashed lines are saddles, dot-dashed lines are totally repelling. On dotted lines the stability is not determined.}
\label{fig:stabFold}
\end{figure}

 We can now address the stability in some cases. In the next result by a {\em totally repelling} fixed point we mean a point that is attracting in reversed time.
 
 \begin{proposition}\label{prop:stabFold}
For $\delta$ close to $\Phi$ and $y=\tau$ the stability of branches created at the fold point with smaller $\tau$ is as follows (see Figure~\ref{fig:stabFold}):
\begin{itemize}
\item
 the branch bifurcating with larger $s\in[0,2\pi)$ consists of attracting fixed points of $\GG$ and the branch bifurcating with lower $s$ created at the same fold point consists of saddles in the following situations:
\begin{itemize}
\item
if $\delta>\Phi$ and $k>1$ near the supercritical fold  in region \ref{item:C};
\item
if $\delta>\Phi$ and $0<k<1$  near the supercritical fold  in region \ref{item:Z};
\end{itemize}
\item
if  $\delta<\Phi$ and $k<1$  near the supercritical  fold on the curve with smaller $\tau$  in region \ref{item:Y}, then
 the branch bifurcating with larger $s\in[0,2\pi)$ consists of totally repelling fixed points of $\GG$ and
  the branch bifurcating with smaller $s$ created at the same fold point consists of saddles;
  \item
if  $\delta<\Phi$ and $k>1$  near the  subcritical  fold on the curve with smaller $\tau$  in region \ref{item:B}, then
  the branch bifurcating with smaller $s\in[0,2\pi)$ consists of totally repelling fixed points of $\GG$
  and   the branch bifurcating with larger $s$ created at the same fold point consists of saddles.
  %
%
%%
%\item
% the branch bifurcating with larger $s\in[0,2\pi)$ consists of totally repelling fixed points of $\GG$ and the branch bifurcating with smaller $s$ created at the same fold point consists of saddles if  $\delta<\Phi$ and $0<k<1$  near the supercritical  fold on the curve with smaller $\tau$  in region \ref{item:Y};
%  %
%  \item
% the branch bifurcating with smaller $s\in[0,2\pi)$ consists of totally repelling fixed points of $\GG$ and the branch bifurcating with larger $s$ created at the same fold point consists of saddles if  $\delta<\Phi$ and $k>1$  near the  subcritical  fold on the curve with smaller $\tau$  in region \ref{item:B}.
%  %
\end{itemize}
\end{proposition}

\begin{proof}
We start by determining the eigenvalues at fold points $\left(y^*,s^*\right)$  of the bifurcation diagram for $\GG$ that occur when $\sin s^*= \pm 1=\varepsilon$ and  $\cos s^*=0$.
At these points the derivative of the map $\GG$ is
\begin{equation}\label{eq:DGfold}
 D\GG \left(y^*,s^*\right)=
 \left(
\begin{array}{ll} 
\delta^2  (y^*)^{\delta^2  -1} +\gamma(\delta^2-\delta)(y^*)^{-p(\delta)}\left( 1+ \varepsilon k\right)&0\\
&\\
\dpt-\frac{1}{K y^*}&1
 \end{array} \right).
\end{equation}
The eigenvalues of $ D\GG{} \left(y^*,s^*\right)$ are $\lambda_1=1$,  as expected at a fold point, and
\begin{eqnarray}
\label{eq:eigenvalue}
\lambda_2(y^*,s^*)&=&\delta^2  (y^*)^{\delta^2  -1} +\gamma(\delta^2-\delta)(y^*)^{-p(\delta)}\left( 1+  \varepsilon k\right)\\
\nonumber&=& (y^*)^{-p(\delta)}\left(\delta^2(y^*)^\delta+\gamma(1+ \varepsilon k)(\delta^2-\delta) \right) >0.
\end{eqnarray}
Therefore, $\lambda_2<1$ if and only if
$\delta^2(y^*)^\delta+\gamma(1+ \varepsilon k)(\delta^2-\delta) < (y^*)^{p(\delta)}$ 
or, equivalently,
$$
\gamma(1+ \varepsilon k)<\dfrac{(y^*)^{p(\delta)}-\delta^2(y^*)^\delta}{\delta^2-\delta}=h_\delta(y^*) .
$$
There are the following two cases to consider:  (I) $\delta>\Phi$ and (II) $\delta<\Phi$.
We now use the properties of $h_\delta$ described in Lemma~\ref{Lemma_aux}.

 \medbreak
 Case (I):   If $\delta>\Phi$,   then $p(\delta)<0$ implying $\dpt\lim_{y^*\to 0+}\lambda_2(y^*,s^*)=0$.
 By assertion \eqref{item:h3} of  Lemma~\ref{Lemma_aux} we have $\dpt\lim_{y\to 0}h_\delta(y)=+\infty$.
Therefore
$h_\delta(y^*)>0$ for all values of $y^*=\tau<\delta^{2/(1-\delta^2)}$ by assertion \eqref{item:h2} of  Lemma~\ref{Lemma_aux}.
Hence, $\lambda_2<1$  corresponds to cases \ref{item:C} and \ref{item:Z} as in Figures~\ref{fig:k>1} and \ref{fig:k<1}, respectively,
with $\varepsilon=1$ and $\gamma(1+k)>1$.
In case  \ref{item:C}  the coordinate $\tau=y^*$ of the fold point  tends to 0 as $\delta$ decreases to $\Phi$ as noted in Subsection~\ref{subsec:k>1} after the proof of Theorem~\ref{th:k>1}.
In case \ref{item:Z} there are two fold points and this analysis holds for the fold point with smaller $\tau$, see Figure~\ref{fig:stabFold}.

  \medbreak
Case (II): If $\delta<\Phi$, then $p(\delta)>0$ implying $\dpt\lim_{y^*\to 0+}\lambda_2(y^*,s^*)=+\infty$.
Therefore at fold points with small $y^*=\tau$ we get $\lambda_2(y^*,s^*)>1$.
For $0<k<1$ fold points with small $\tau$ as $\delta$ increases to $\Phi$ occur at the supercritical fold with $\gamma(1-k)<1$  and $\varepsilon=-1$ in region \ref{item:Y} (see Figure~\ref{fig:Foldsk<1})
whereas for $k>1$ they take place at the subcritical
fold with $\gamma(1+k)<1$  and  $\varepsilon=1$  in region \ref{item:B} (Figure~\ref{fig:Foldsk>1}).

 \medbreak
Finally, we  estimate  the eigenvalues  of $D\GG(y,s)$ at points $(y,s)$ in the branches bifurcating from the fold points described above.
Since $s$ is close to $s^*=(2-\varepsilon)\pi/2$, with $\varepsilon=\pm 1$, then 
$$
\cos s= -\varepsilon(s-s^*)+\mathcal{O}\left((s-s^*)^3\right)
\quad\mbox{and}\quad 
\sin s=\varepsilon+\mathcal{O}\left((s-s^*)^2\right) ,
$$
where $\mathcal{O}$ is the standard Landau notation.
Hence we get the following expressions for the determinant and the trace of the matrix  $D\GG(y,s)$:
$$
\begin{array}{lcl}
\det D\GG(y,s)&=&\delta^2  y^{\delta^2  -1} +\gamma(\delta^2-\delta)y^{-p(\delta)}\left( 1+\varepsilon k\right)
-\dfrac{\varepsilon\gamma k}{K}y^{-p(\delta)}(s-s^*)+\mathcal{O}\left((s-s^*)^2\right)\\
{}&\approx&\lambda_2(y,s)-\dfrac{\varepsilon\gamma k}{K}y^{-p(\delta)}(s-s^*)+\mathcal{O}\left((s-s^*)^2\right)
\end{array}
$$
 and 
 $$
 \tr D\GG(y,s)\approx 1+ \lambda_2(y^*,s^*)+\mathcal{O}\left((s-s^*)^2\right).
 $$
  When $s$ moves away from $s^*$ the trace of  $D\GG(y,s)$ remains roughly the same.

  \medbreak
   If $\delta>\Phi$ with  $\varepsilon=1$, since $-\dfrac{\varepsilon\gamma k}{K}y^{-p(\delta)}<0$, then   the determinant $\det D\GG(y,s)$ increases for small $s<s^*$ and decreases for small $s>s^*$, so the fixed point that bifurcates with larger $s$ is attracting.
 
 \medbreak 
  If $1<\delta<\Phi$  and   $k>1$ with $\varepsilon=1$, then again  the determinant $\det D\GG(y,s)$ decreases for small $s>s^*$ and increases for small $s<s^*$. 
 The fixed point bifurcating with smaller $s$ is totally repelling while the fixed point with larger $s$ is a saddle.

 \medbreak 
 If $1<\delta<\Phi$  and   $0<k<1$ with $\varepsilon=-1$, then near  $s^*$ the  determinant increases with $s$, hence the totally repelling fixed point is the one with larger $s$ while the fixed point with smaller $s$ is a saddle.
  Thus we obtain the stability conditions in the statement. 
\end{proof}

\subsection{A  bifurcation of codimension 2}\label{subsec:cod2}
 In Subsection \ref{ss: SN} we have assigned the stability of the   fixed points
 emerging at the two branches of the saddle-node bifurcations in regions  \ref{item:B}, \ref{item:C},  \ref{item:Y} and  \ref{item:Z}. Here we show that these codimension 1 bifurcations are part of a generic discrete-time Bogdanov--Takens bifurcation  in the parameter space $(\tau,\delta,\gamma)$.
This bifurcation takes place when the derivative   of $\GG$  has 1 as a double eigenvalue and is not the identity map.
The next result shows that this happens when
 two fold points come together at the maximum $\tau_m$ of $\FF_\delta(\tau)$,
where $\FF_\delta(\tau_m)=M_\FF(\delta)$, as defined in \eqref{item:Fdelta<phi} of Lemma~\ref{lem:Fdelta}.
These are precisely the transitions at Morse critical points of $g(\tau,s)$ namely 
\ref{item:W}$\leftrightarrow$\ref{item:X},
\ref{item:X}$\leftrightarrow$\ref{item:Y},
\ref{item:A}$\leftrightarrow$\ref{item:B} and
\ref{item:a}$\leftrightarrow$\ref{item:b}.

\begin{theorem}
\label{teoremaEstabilidade1} 
For each $\delta$ with $1<\delta<\Phi$ and any $k>0$, if $\gamma = M_\FF(\delta)/(1+  k)$, 
then there is a  fold point  $\left(y^*,s^*\right)$ for the map $g$  defined in \eqref{eq:g0} such that
the derivative $D\GG\left(y^*,s^*\right)$  of $\GG$
has 1 as a double eigenvalue and is not the identity,
with $y^*=\tau_m(\delta)$ and $\sin s^*=1$.

Moreover, if $0<k<1$, then the same is true for $\gamma = M_\FF(\delta)/(1-  k)$ with $\sin s^*=-1$.
\end{theorem}

\begin{proof} 
A fold point $\left(y^*,s^*\right)$ for $g$ must satisfy $\sin s^*=\pm 1=\varepsilon $, as discussed in the proof of Proposition~\ref{prop:stabFold}.
First we show that under these conditions $y^*=\tau_m$ if and only if $\gamma = M_\FF(\delta)/(1+ \varepsilon k)$.
For $y^*=\tau_m$ and $\sin s^*=\pm 1=\varepsilon$,   by \eqref{eq:Fdelta}, we have 
$g\left(\tau_m,s^*\right)=0$ if and only if $\gamma(1+ \varepsilon  k)=\FF_\delta(\tau_m)$.
From Lemma~\ref{lem:Fdelta} it follows that $\FF_\delta(y^*)=M_\FF(\delta)$.

Since $\sin s^*=\pm 1$, then $\cos s^*=0$ and hence 
$\dfrac{\partial g}{\partial s}\left(\tau_m,s^*\right)=0$.
As in the proof of Proposition~\ref{prop:stabFold},  the matrix  $D\GG\left(y^*,s^*\right)$ is triangular and 
its  eigenvalues  are  $\lambda_1=1$ and $\lambda_2$ 
 where $k$ is replaced by $\varepsilon k$ in the expression \eqref{eq:eigenvalue}.
 
Since $\left(y^*,s^*\right)$ is a fixed point of $\GG$, then equation \eqref{eq:Fdelta} holds and hence
$\gamma(1+ \varepsilon k)=(y^*)^{p(\delta)}-(y^*)^\delta =\FF_\delta(y^*)$.
Replacing this value of $\gamma(1+  \varepsilon k) $ in the value of $\lambda_2$ in \eqref{eq:eigenvalue} and using 
the expression for $y^*=\tau_m$ in assertion \eqref{item:Fdelta<phi} of Lemma~\ref{lem:Fdelta} we get
\begin{equation}\label{eqEigenvs}
\lambda_2(y^*,s^*)=\delta  (y^*)^{\delta^2  -1} +\delta^2-\delta =  \delta \left(\dfrac{-\delta^2+\delta+1}{\delta} \right)^{(\delta^2-1)/(\delta^2-1)}+\delta^2-\delta =1.
\end{equation}
The result on $D\GG$ follows since the matrix \eqref{eq:DGfold} is not diagonal.
 \end{proof}
 
 For  $0<k<1$, the parameters $\gamma$ and $\delta$ for which the conditions of Theorem~\ref{teoremaEstabilidade1}  hold when  $\varepsilon=-1$ correspond to the transition between cases \ref{item:W} and \ref{item:X}  of Theorem~\ref{th:k<1}. 
 When $\varepsilon=+1$ they lie on the transition between cases \ref{item:X} and \ref{item:Y}, as indicated in Figure~\ref{fig:k<1}.
 When $k>1$ only the case $\varepsilon=+1$ occurs and it corresponds to the  transition between cases \ref{item:A} and \ref{item:B}  of Theorem~\ref{th:k>1} shown in Figure~\ref{fig:k>1}.
 For $k=1$ this is the transition between cases \ref{item:a} and \ref{item:b}  of Proposition~\ref{th:k=1}.

\begin{figure}[hht]
\includegraphics[width=109mm]{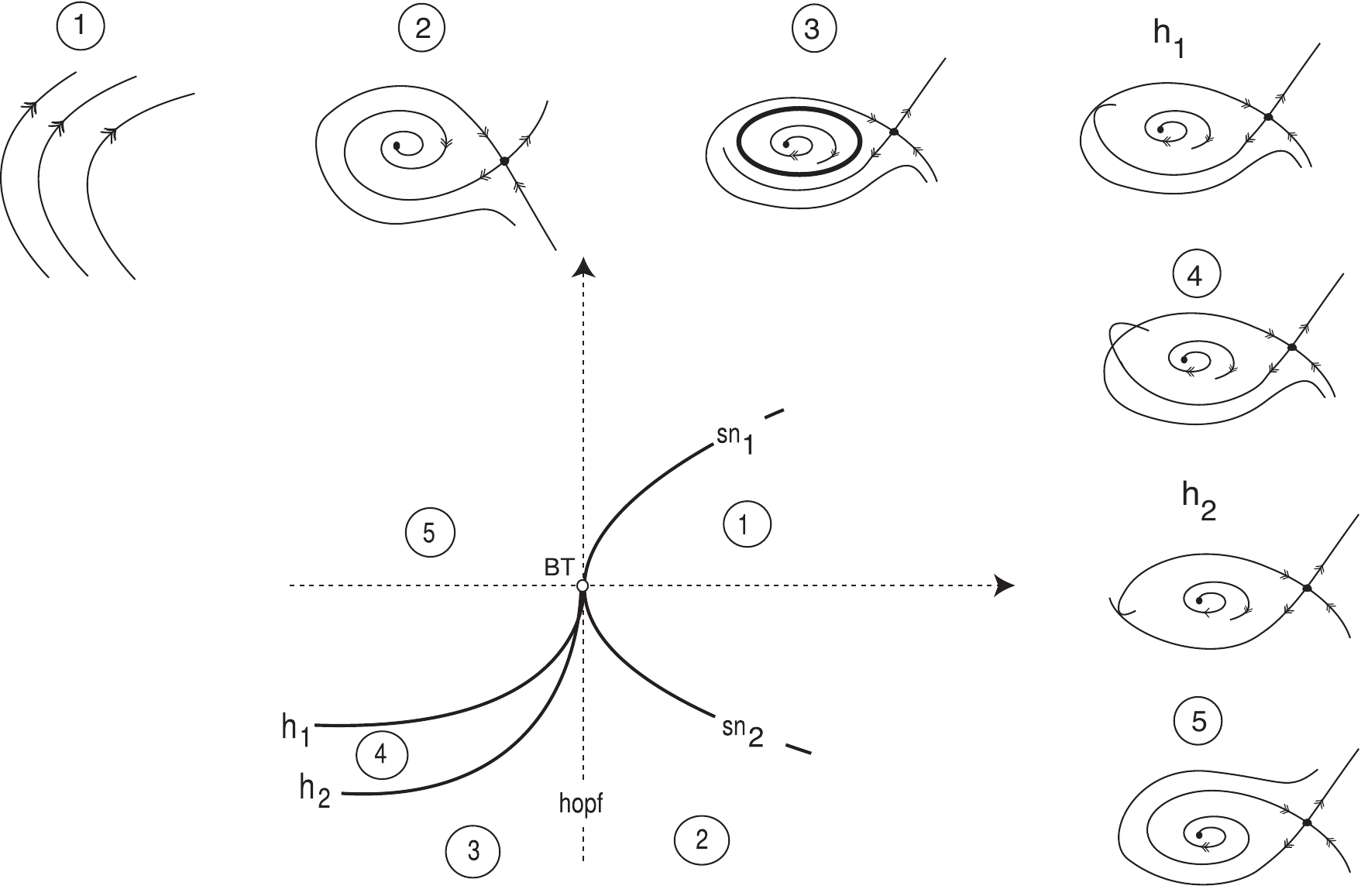}\\
 \caption{\small Dynamics for a generic discrete-time Bogdanov--Takens bifurcation. 
Curves:  $sn_1$ and $sn_2$ saddle-node bifurcations; $hopf$ Hopf bifurcation; $h_1$ and $h_2$ homoclinic tangencies associated to a dissipative fixed point. Regions: $(1)$   no recurrent dynamics; $(2)$   the unstable manifold of the saddle is connected to the stable manifold of the focus; $(3)$ a saddle and a stable periodic orbit; $(4)$   horseshoe dynamics; $(5)$  the unstable manifold of the focus intersects the stable manifold of the saddle.}
\label{fig:BT4b}
\end{figure}

Fold points $\left(y^*,s^*\right)$ satisfy $\sin s^*=\pm 1=\varepsilon$.
For given fixed $\tau\in(0,1]$, $\gamma>0$ and $\delta>1$ there is a fold point if $y^*=\tau$ 
with $\FF_\delta(\tau)=\gamma(1+\varepsilon k)$.
These are shown in Figures~\ref{fig:Foldsk>1} and \ref{fig:Foldsk<1} as curves in the $(\tau,\delta)$ plane for fixed $\gamma$.
As $\gamma$ varies this forms a surface in the $(\tau,\delta,\gamma)$ parameter space, a set of codimension 1.
Theorem~\ref{teoremaEstabilidade1} implies that this surface of fold points intersects the set $(\tau,\delta,\gamma)=(\tau_m(\delta),\delta,\gamma)$ at two curves (a codimension 2 set), where
 $$
\tau=\tau_m(\delta)\qquad\gamma = M_\FF(\delta)/(1\pm k)\qquad 1<\delta<\Phi 
 \qquad 0<k<1 .
$$
Theorem~\ref{teoremaEstabilidade1} also shows that, for $(\tau,\delta,\gamma)$ on these curves and at the fold point $\left(y^*,s^*\right)$ where it holds, the derivative $D\GG\left(y^*,s^*\right)$ 
has 1 as a double eigenvalue and is not the identity.
At these points we expect to find a discrete-time Bogdanov--Takens bifurcation \cite{BRS,Yagasaki} -- cf. Figure \ref{fig:BT4b}.
This is a codimension 2 bifurcation occurring at points where 1 is a double eigenvalue,
the derivative is not the identity and where the map   satisfies some nondegeneracy conditions on the nonlinear part.
Instead of verifying these additional conditions we check that the linear conditions for nearby local bifurcations arise in a form consistent with the versal unfolding of the discrete-time Bogdanov--Takens bifurcation.
For parameter values around this bifurcation, by the results of \cite{BRS,Yagasaki}, we expect to find the following sets in the $(\tau,\delta,\gamma)$ parameter space:
\begin{enumerate}
\renewcommand{\theenumi}{(\alph{enumi})}
\renewcommand{\labelenumi}{{\theenumi}}
\item\label{saddlenodes}  a surface of saddle-node bifurcations;
\item\label{Hopf}   a surface of  Hopf bifurcations;
\item\label{homoclinic}   two surfaces of homoclinic tangencies.
\end{enumerate}

The saddle-node bifurcations \ref{saddlenodes} are the fold points that we have already identified, given by $g(s^*,\tau)=0$ with $\sin s^*=\pm 1$. 
The corresponding surface has two components defined by the two conditions 
$\FF_\delta(\tau)=\gamma(1\pm k)$.

The surface defined by the linear conditions for Hopf bifurcations \ref{Hopf} near these points is treated in 
Subsection~\ref{subsec:Hopf} below.
 Generically a $\GG$-invariant circle bifurcates at this surface.

The two surfaces \ref{homoclinic}  correspond to bifurcations at which the stable and unstable manifolds of a saddle point are tangent.
In the region between these surfaces 
there is a transverse intersection of the stable and the unstable manifolds of the saddle.
These intersections are repeated at an orbit that accumulates on the saddle in forward and backward times.
Around the transverse intersection of the surfaces horseshoe dynamics occurs. 
The two surfaces are tangent at the points of codimension two bifurcation and they delimit a region that gets exponentially small as it approaches the  curve of Bogdanov--Takens bifurcation.
The invariant manifolds intersect transversely within the parameter region between these two surfaces and do not intersect outside it.  
This configuration implies that in this region the dynamics of $\GG$  exhibits behaviour
equivalent to Smale's horseshoe and Newhouse phenomena.

\subsection{Hopf bifurcation}\label{subsec:Hopf}
In this section we address statement \ref{Hopf} above that the parameter space $(\tau,\delta,\gamma)$ 
contains a surface of Hopf bifurcation points that intersects any neighbourhood of  the points where $D\GG\left(y^*,s^*\right)$ 
has 1 as a double eigenvalue.
Hopf bifurcation is a codimension 1 bifurcation occurring at points where the derivative $D\GG(\tau,s)$ has non-real eigenvalues of absolute value 1 and where the nonlinear part of $\GG$  also satisfies   nondegeneracy conditions.
In Proposition~\ref{prop:Hopf} below we show that  the condition on the eigenvalues holds on a surface that contains the discrete Bogdanov--Takens points in its closure.
In Corollary~\ref{corol:Hopf} we  describe the set of parameters where the condition  on the eigenvalues is satisfied.

We start with some properties of the trace and determinant of $D\GG(\tau,s)$ that will be used in the proof of Proposition~\ref{prop:Hopf}.

\begin{lemma}\label{lema:detDG}
Let $(\tau,s(\tau))\in\CC$ be a branch in the bifurcation diagram of $g$ (i.e., $g(\tau,s(\tau))\equiv 0$).
Then $\det D\GG(\tau,s(\tau))=1$ if and only if  $f(\tau,\delta,\gamma)=0$ where
\begin{equation}\label{eq:1detDG=1}
f(\tau,\delta,\gamma)=\gamma k\left( \frac{ 1}{K} \cos s(\tau)+(\delta^2-\delta)\sin s(\tau)\right)-
\tau^{p(\delta)} +\delta^2   \tau^{\delta}+\gamma(\delta^2-\delta)
\end{equation}
 and
$p(\delta)$ is the polynomial defined in \eqref{eq:pdelta}.
In this case
$$
\tr D\GG(\tau,s(\tau))=2- \dfrac{ \gamma k}{K} \tau^{-p(\delta)} \cos s(\tau) .
$$

\end{lemma}
\begin{proof}
 Using \eqref{eq:DG} for $y=\tau$ we get 
$$
\det D\GG(\tau,s(\tau))=
\delta^2   \tau^{\delta^2  -1} +\gamma(\delta^2-\delta)\tau^{-p(\delta)}\left( 1+k\sin s(\tau)\right)+
 \frac{ \gamma k}{K} \tau^{-p(\delta)} \cos s(\tau) .
$$
Therefore, $\det D\GG(\tau,s(\tau))=1$
if and only if
\begin{equation}\label{eq:detDG=1}
 \frac{ \gamma k}{K} \tau^{-p(\delta)} \cos s=
1-\left(\delta^2   \tau^{\delta^2  -1} +\gamma(\delta^2-\delta)\tau^{-p(\delta)}\left( 1+k\sin s\right)\right).
\end{equation}
Equivalently, multiplying by $\tau^{p(\delta)}$ and taking into account that $\tau^{p(\delta)}\tau^{\delta^2-1}=\tau^\delta$, we get:
$$
\tau^{p(\delta)} -\delta^2   \tau^{\delta}+\gamma(\delta^2-\delta)=
\gamma k\left( \frac{ 1}{K} \cos s+(\delta^2-\delta)\sin s
\right)
$$
i.e., $f(\tau,\delta,\gamma)=0$ as stated.

Again, using \eqref{eq:DG} for $y=\tau$, we get 
$$
\tr D\GG(\tau,s(\tau))=1+\delta^2   y^{\delta^2  -1} +\gamma(\delta^2-\delta)y^{-p(\delta)}\left( 1+k\sin s(\tau)\right)
$$
and if $\det D\GG(\tau,s(\tau))=1$, then from \eqref{eq:detDG=1} it follows that
$$
\delta^2   \tau^{\delta^2  -1} +\gamma(\delta^2-\delta)\tau^{-p(\delta)}\left( 1+k\sin s(\tau)\right)=
1- \frac{ \gamma k}{K} \tau^{-p(\delta)} \cos s(\tau),
$$
yielding the expression for $\tr D\GG(\tau,s(\tau))$ in the statement.
\end{proof}

The map $s(\tau)$ depends explicitly on $\tau$ which in turn depends on $\gamma$. 
We  use  this dependence in the next result.

\begin{lemma}\label{lema:dfdgama}
For any given $\gamma_0>0$ and $\delta_0>1$, if $(\tau_0,s(\tau_0))$ is a fold point  in the branch $(\tau,s(\tau))\in\CC$ of the bifurcation diagram of $g$ with $\sin s(\tau_0)=\varepsilon=\pm 1$, then
$$
\dfrac{\partial f}{\partial \gamma}(\tau_0,\delta_0,\gamma_0)=
\left(1+\varepsilon k\right)({\delta_0}^2-\delta_0)\ne 0 .
$$
\end{lemma}
\begin{proof}

Since $\sin s(\tau_0)=\varepsilon=\pm 1$, then $\cos s(\tau_0)=0$ and 
$$
\dfrac{\partial f}{\partial \gamma}(\tau,\delta,\gamma)
=k\left( \frac{ 1}{K} \cos s(\tau)+(\delta^2-\delta)\sin s(\tau)\right)
+(\delta^2-\delta)+
\gamma k\left( \frac{ -1}{K} \sin s(\tau)+(\delta^2-\delta)\cos s(\tau)\right)
\dfrac{\partial s(\tau)}{\partial \gamma} ,
$$
implying that
\begin{equation}\label{eq:dfdgamma}
\dfrac{\partial f}{\partial \gamma}(\tau_0,\delta_0,\gamma_0)=
\left(1+\varepsilon k\right)({\delta_0}^2-\delta_0) -
\gamma_0 k \frac{ \varepsilon}{K}\dfrac{\partial s(\tau)}{\partial \gamma}.
\end{equation}
In order to estimate $\dfrac{\partial s(\tau)}{\partial \gamma}$ we use $g(\tau,s(\tau))\equiv 0$ from \eqref{eq:g0}.
Therefore,
$$
\dfrac{\partial g}{\partial \gamma}(\tau,s(\tau))=
\tau^{\delta^2-\delta}\left(1+k\sin s(\tau)\right)+\gamma k\tau^{\delta^2-\delta}\cos s(\tau)\dfrac{\partial s(\tau)}{\partial \gamma}\equiv 0
$$
and
$$
\dfrac{\partial^2 g}{\partial \gamma^2}(\tau,s(\tau))=
k\tau^{\delta^2-\delta}\cos s(\tau)\dfrac{\partial s(\tau)}{\partial \gamma}+
\gamma k\tau^{\delta^2-\delta}\left(-\sin s(\tau)  \dfrac{\partial s(\tau)}{\partial \gamma}+\cos s(\tau)\dfrac{\partial^2 s(\tau)}{\partial \gamma^2}\right)\dfrac{\partial s(\tau)}{\partial \gamma}\equiv 0 ,
$$
hence
$$
\left.
\dfrac{\partial^2 g}{\partial \gamma^2}(\tau_0,s_0)
\right|_{(\delta_0,\gamma_0)}=
-\varepsilon\gamma_0 k{\tau_0}^{{\delta_0}^2-\delta_0}  \left(\dfrac{\partial s(\tau_0)}{\partial \gamma}\right)^2=0
\quad\Rightarrow\quad
\left.\dfrac{\partial s(\tau_0)}{\partial \gamma}\right|_{(\delta_0,\gamma_0)}=0 .
$$
Substituting into \eqref{eq:dfdgamma} we obtain, for $(\tau_0,\delta_0,\gamma_0)$ satisfying the conditions of Theorem~\ref{teoremaEstabilidade1},  that
$$
\dfrac{\partial f}{\partial \gamma}(\tau_0,\delta_0,\gamma_0)=
\left(1+\varepsilon k\right)({\delta_0}^2-\delta_0)\ne 0.
$$
Note that the case $\varepsilon=-1$ does not occur for $k=1$.
\end{proof}

For the next result recall that the expressions of $g(\tau,s)$ and of $\GG(\tau,s)$ also depend on the parameters  $\gamma$ and $\delta$.

\begin{proposition}\label{prop:Hopf}
 If $\sin s^*=\pm 1$ with $s^*\in[0,2\pi)$, if $(\tau^\star,\delta^\star,\gamma^\star)$ satisfies the conditions of Theorem~\ref{teoremaEstabilidade1} and $k>0$, then $(\tau^\star,\delta^\star,\gamma^\star,s^*)$ is an accumulation point of
 % $H\subset [0,1]\times (1,+\infty)\times\RR_+\times[0,2\pi)$
$$
H=\left\{(\tau,\delta,\gamma,s):
\ g(\tau,s)=0\mbox{ and for } y=\tau \ 
D\GG(y,s)\mbox{ has eigenvalues } \ee^{\pm i\theta}\mbox{ with } \theta \in(0,\pi)
\right\}.
$$
Moreover, if ${\mathbf p}:\RR^4\seta\RR^3$ is the projection ${\mathbf p}(\tau,\delta,\gamma,s)=(\tau,\delta,\gamma)$,  then ${\mathbf p}(H)$
contains a smooth surface in $\RR^3$.
\end{proposition}

\begin{proof}
For any $\delta>\delta^\star$ with small  $\delta-\delta^\star$ and any  $\gamma$ such that $\gamma^\star-\gamma>0$ is small the bifurcation diagram for $g(\tau,s)=0$ has a fold point near $(\tau^\star,s^*)$ with  $\cos s>0$ on one of the branches containing the fold.
From now on we assume $(\tau,s)=(\tau,s(\tau))$ to be on this branch.
The idea of the proof is to show that 
given $(\tau,s(\tau))$ close to the fold point on this branch
we may find 
a smooth function 
$\gamma=\gamma(\tau,\delta)$ close to $\gamma^\star$ such that 
$(\tau,\delta,\gamma(\tau,\delta),s(\tau))\in H$.

The derivative $D\GG(\tau,s)$ has eigenvalues $ \ee^{\pm i\theta}$ with $\theta \in(0,\pi)$ if and only if $\det D\GG(\tau,s)=1$ and  $-2<\tr D\GG(\tau,s)<2$.
At $(\tau^\star,\delta^\star,\gamma^\star)$ we have, from Theorem~\ref{teoremaEstabilidade1}:
$$
\det D\GG(\tau^\star,s^*)=1 \quad \text{and} \quad \tr D\GG(\tau^\star,s^*)=2.
$$
By Lemma~\ref{lema:detDG} we have $\det D\GG(\tau,s(\tau))=1$ if and only if  $f(\tau,\delta,\gamma)=0$.
Since  we have established in Lemma~\ref{lema:dfdgama} that $\dfrac{\partial f}{\partial \gamma}(\tau^\star,\delta^\star,\gamma^\star)\ne 0$, therefore,
by the implicit function theorem, we get $\det D\GG(\tau,s(\tau))=1$ for some function $\gamma=\gamma(\tau,\delta)$
with $\gamma(\tau^\star,\delta^\star)=\gamma^\star$.

It remains to check that at these points we have  $-2<\tr D\GG(\tau,s(\tau))<2$.
Applying Lemma~\ref{lema:detDG}
to the branch of $g(\tau,s)=0$ where  $\cos s>0$ that was chosen in the beginning we have
$$
\tr D\GG(\tau,s(\tau))=2- \dfrac{ \gamma k}{K} \tau^{-p(\delta)} \cos s(\tau)<2.
$$

Since $\tr D\GG(\tau^\star,s^*)=2$ at $(\tau^\star,\delta^\star,\gamma^\star)$, then 
 $\tr D\GG(\tau,s(\tau))>-2$ for $(\tau,\delta,\gamma)$ close to this point.

Finally, we have established that given $\gamma=\gamma(\tau,\delta)$ there exists $s=s(\tau)$ such that 
$(\tau,\delta,\gamma(\tau,\delta),s(\tau))\in H$.
It follows that the smooth surface defined by $\{(\tau,\delta,\gamma(\tau,\delta))\}$ is contained in ${\mathbf p}(H)$.
\end{proof}

We had claimed in \ref{Hopf} that there was a surface in  the parameter space $(\tau,\delta,\gamma)$ of Hopf bifurcation points. 
The set $H$ lies in $(\tau,\delta,\gamma,s)$ space and we have shown that the graph of $\gamma(\tau,\delta)$ is contained in its projection into $(\tau,\delta,\gamma)$ space.
Therefore, it should be the surface in the claim but we cannot guarantee that   all the points of ${\mathbf p}(H)$ lie on this graph.

From the proofs of Proposition~\ref{prop:stabFold}, of Theorem~\ref{teoremaEstabilidade1} and of Proposition~\ref{prop:Hopf} we can locate in parameter space  the sites of Hopf bifurcation points  where an invariant circle emerges   from a fixed point.
This is the content of the next result.
A  Hopf bifurcation at a point $(\tau_H,s_H)$ can be {\em subcritical} if the invariant circle exists for $\tau<\tau_H$, or {\em supercritical} if the  circle exists for $\tau>\tau_H$.

\begin{corollary}\label{corol:Hopf}
A  Hopf bifurcation  for the map $g$  defined in \eqref{eq:g0}  creating an  invariant circle  is possible for all $\delta$ with $1<\delta<\Phi$ for the following sets of the parameter $\gamma$:
  \begin{itemize}
  \item
  for $k>1$  and
  $0<\gamma<\gamma_+(\delta)$ close to $\gamma_+(\delta)$ (region \ref{item:B} of Theorem~\ref{th:k>1}) on the branch with a subcritical fold at $\sin s^*=1$ with $s<s^*$;
   \item
  for $0<k<1$ and
  $\gamma_+(\delta)<\gamma<\gamma_-(\delta)$ close to either $\gamma_-(\delta)$  or $\gamma_+(\delta)$ (region \ref{item:X} of Theorem~\ref{th:k<1}) on the branch with a supercritical fold at $\sin s^*=1$ with $s>s^*$; 
    \item
  for $0<k<1$ and
   $0<\gamma<\gamma_+(\delta)$ close to $\gamma_+(\delta)$ (region \ref{item:Y} of Theorem~\ref{th:k<1}) on the branch with a supercritical fold at $\sin s^*=-1$ with $s>s^*$;   ,
  \end{itemize}
  where $\gamma_+(\delta)=M_{\FF}(\delta)/(1+k)$ and $\gamma_-(\delta)=M_{\FF}(\delta)/(1-k)$ and
$M_\FF(\delta)>0$ is defined  in assertion \eqref{item:Fdelta<phi} of Lemma~\ref{lem:Fdelta}.

The invariant circle created in  region \ref{item:B} will be attracting if the  Hopf bifurcation  is supercritical. 
The circle in region \ref{item:Y} will be attracting if the bifurcation  is subcritical. 
\end{corollary}

\begin{proof}
The proof of Proposition~\ref{prop:Hopf} consists of finding  a condition on the parameters $\tau$, $\delta$ and $\gamma$ under which  the derivative $D\GG(y,s)$ with $y=\tau$ at a fixed point of $\GG(y,s)$ has a non real eigenvalue at the unit circle.
The condition is that $f(\tau,\delta,\gamma(\tau,\delta))=0$ where $f$ is defined by \eqref{eq:1detDG=1} of Lemma~\ref{lema:detDG}.
The function $\gamma(\tau,\delta)$ is obtained applying the implicit function theorem at a point $(\tau^*,\delta^*,\gamma^*)$ where Theorem~\ref{teoremaEstabilidade1} holds.
In particular, $1<\delta^*<\Phi$ and either $\gamma^*=\gamma_+(\delta)$ or $\gamma^*=\gamma_-(\delta)$ in the notation of the present statement.
Calculations are done at a fold point $(\tau,s^*)$ where $\sin s^*=\pm 1$ if $0<k<1$ and $\sin s^*=+ 1$ if $k>1$.
This means that $(\delta,\gamma)$ must either be in region \ref{item:B} of Theorem~\ref{th:k>1} or
in one of the regions \ref{item:X} or \ref{item:Y} of Theorem~\ref{th:k<1}, since the bifurcation diagrams corresponding to region  \ref{item:A}  do not contain fold points and those in region \ref{item:W} do not contain fixed points.

Proposition~\ref{prop:stabFold} established that near the fold points in regions \ref{item:B} and \ref{item:Y} the branches in the bifurcation diagrams concerned in the proof of Proposition~\ref{prop:Hopf} are totally repelling.
At the Hopf bifurcation the branch becomes attracting.
The invariant circle will be attracting if it bifurcates towards the fold point.
\end{proof}

\section{Frequency locking}\label{sec:freq}

In this section we  describe the interpretation for the original problem \eqref{general} of results in the previous sections. 
These arise as frequency locked solutions, a typical phenomenon in periodically perturbed differential equations,
 as reported numerically by Tsai and Dawes \cite{TD1} in the context of the May--Leonard system.
 %, for instance.
We describe the bifurcation of these solutions when the frequency $\omega$ of the forcing term is varied.

\begin{definition}\label{def:frequencyLock}
A periodic solution of \eqref{general} is said to be \emph{frequency locked} if its period is an integer multiple of the period of the external forcing.
\end{definition}

Returning to the original problem  we look for {\em frequency locked} solutions of \eqref{general}: periodic solutions  whose period is an integer multiple of $\pi/\omega$, the period of the forcing term.

A point $(y,s)\in\CC$ corresponds to a periodic solution of  \eqref{general} with period $P$
if $G(y,s)=(y,s+P)$, i.e., 
 if $(y,s)$ is a fixed point of $$\GG(y,s)=G(y,s)+(0,\ln\tau/K) \quad \text{with}\quad P=-\ln\tau/K.$$
For each forcing frequency $\omega>0$ there are frequency locked solutions when $\tau=\ee^{-nK\pi/\omega}$, $n\in\NN$, where $1:n$ is the {\em frequency ratio} of the solution.
Hence $1:n$ frequency locked solutions occur when $\omega=-nK\pi/\ln \tau$.
Therefore, there are $1:1$ frequency locked solutions for some $\omega=\omega_0>0$ if and only if  there are $1:n$ frequency locked solutions for  $\omega=\omega_0/n$, $n\in\NN$.

Consider fixed values of $\delta>\Phi$ of $\gamma>0$  and of $0<k< 1$.
Then, from Theorem~\ref{th:k<1}, there are two solutions to $g(\tau,s)=0$ for each $\tau$ in some open interval $(\tau_1,\tau_2)$ with $0<\tau_1<\tau_2<1$.
Therefore there are $1:1$ frequency locked solutions of \eqref{general} for each $\omega$ between $\omega_1=-K\pi/{\ln\tau_1}$ and $\omega_2=-K\pi/{\ln\tau_2}$.
Also, if $(\tau,s)$ is a point in the bifurcation diagram of $g(\tau,s)=0$, then there is a $1:n$ 
frequency locked solution of \eqref{general} for $ \omega=-nK\pi/\ln \tau$ with $ n\in \NN$.
 The next result follows from applying the same reasoning to the other regions discussed in Theorem~\ref{th:k>1}.

\begin{corollary}\label{cor:k<1}
If $0<k<1$, then there are  real numbers $0<\omega_1<\omega_2<\omega_3<\omega_4$ such that
there are two $1:1$ frequency locked periodic solutions of \eqref{general} under any of 
the following conditions:
\begin{enumerate}
\renewcommand{\theenumi}{(\Alph{enumi})}
\renewcommand{\labelenumi}{{\theenumi}}
\setcounter{enumi}{23}
\item
for $\omega\in(\omega_1,\omega_2)$ if $\delta<\Phi$ and $ \gamma_+(\delta)=M_\FF(\delta)/(1+k)<\gamma< M_\FF(\delta)/(1-k) =\gamma_-(\delta)$;
\item
for $\omega\in(\omega_1,\omega_2)\cup(\omega_3,\omega_4)$ if $\delta<\Phi$ and $0<\gamma<M_\FF(\delta)/(1+k) =\gamma_+(\delta)$;
\item
for $\omega\in(\omega_1,\omega_2)$ if $\delta>\Phi$.
\end{enumerate}
Moreover, for $n\in\NN$ there are $1:n$ 
frequency locked periodic solutions for some value of $\omega$ if and only if there are 
$1:1$ frequency locked solutions for $n\omega$.

{\rm \ref{item:W}} If $\delta<\Phi$ and $\gamma_-(\delta)=M_\FF(\delta)/(1-k)<\gamma$, then there are no frequency locked periodic solutions of \eqref{general}.
\end{corollary}

Similarly, from Theorem~\ref{th:k<1} we get:

\begin{corollary}\label{cor:k>1}
If $k>1$, then  there are  real numbers $0<\omega_1<\omega_2$ such that
 there are two $1:1$ frequency locked periodic solutions of \eqref{general}  under any of
 the following conditions:
\begin{enumerate}
\renewcommand{\theenumi}{(\Alph{enumi})}
\renewcommand{\labelenumi}{{\theenumi}}
\item
for all $\omega\in\RR_+$ if $\delta<\Phi$ and $\gamma> M_\FF(\delta)/(1+k) =\gamma_+(\delta)$;
\item
for $\omega\in(\omega_1,\omega_2)$ if $\delta<\Phi$ and $0<\gamma<M_\FF(\delta)/(1+k) =\gamma_+(\delta)$;
\item
for all $\omega>\omega_1$ if $\delta>\Phi$.
\end{enumerate}
Moreover, for $n\in\NN$ there are $1:n$ 
frequency locked periodic solutions for some value of $\omega$ if and only if there are 
$1:1$ frequency locked solutions for $n\omega$.
\end{corollary}

Consider fixed values of $\delta>1$, $\gamma>0$ and $k>0$  in one of the regions described by one of Corollaries~\ref{cor:k<1} or \ref{cor:k>1} and let $(\tau_0,s_0)$ be  any point in the bifurcation diagram of $g(\tau,s)$.
Then there is a $1:1$ frequency locked solution for $\omega_0=-K\pi/\ln\tau_0$ with period $P=-\ln\tau_0/K$.
When $\omega$ grows from $\omega_0$ there appear $1:n$ frequency locked solutions at $\omega=n\omega_0$ with period $P=-n\ln\tau_0/K=n\pi/\omega_0$.
The period $\pi/\omega$ of the forcing is reduced and the periods  of the bifurcating  frequency locked solutions increase, since they have a larger frequency ratio.
Hence, reducing the forcing period  leads to the bifurcation  of solutions with arbitrarily large periods.

\medbreak

Definition~\ref{def:frequencyLock} may be extended to invariant sets for the differential equation  \eqref{general}. 
A flow-invariant set ${\mathcal L}({\mathcal S})$ for \eqref{general} is the suspension  of a $\GG$--invariant set ${\mathcal S}$ on $\CC$.
When $-\ln \tau/K$ is an integer multiple of the period $\pi/\omega$ of the forcing term, then solutions lying on the set ${\mathcal L}({\mathcal S})$ are also said to be frequency locked.

The next result summarises the consequences of the results in 
Section~\ref{sec:stability} for the case  $\delta\in\left(1,\Phi\right)$. 

\begin{proposition}\label{prop:deltaSmall}
For any $\delta\in\left(1,\Phi\right)$ let $\gamma_+(\delta)=M_{\FF}(\delta)/(1+k)$ and $\gamma_-(\delta)=M_{\FF}(\delta)/(1-k)$.
If either $k>1$ or $0<k<1$ and $\gamma<\gamma_-(\delta)$, then there are frequency locked solutions of  \eqref{general}.

Moreover, for the values of the parameter $\gamma$ listed below there is
a frequency locked invariant torus arising  at a Hopf bifurcation from a frequency locked periodic solution of  \eqref{general} that in turn bifurcates into a frequency locked suspended horseshoe in the small region in parameter space delimited by the two homoclinic tangencies shown in Figure~\ref{fig:BT4b}.
The values of $\gamma$ where the invariant torus exists are:
  \begin{itemize}
  \item
for $k>1$  and
  $0<\gamma<\gamma_+(\delta)$ close to $\gamma_+(\delta)$ (region \ref{item:B} of Theorem~\ref{th:k>1});
   \item
for $0<k<1$ and  $\gamma_+(\delta)<\gamma<\gamma_-(\delta)$ close to either $\gamma_-(\delta)$  or $\gamma_+(\delta)$ (region \ref{item:X} of Theorem~\ref{th:k<1});
    \item
for $0<k<1$ and   $0<\gamma<\gamma_+(\delta)$ close to $\gamma_+(\delta)$ (region \ref{item:Y} of Theorem~\ref{th:k<1}).
  \end{itemize}
\end{proposition}

The next result contains the consequences of Theorems~\ref{th:k>1}--\ref{th:k=1} and Proposition~\ref{prop:stabFold} for the case $\delta>\Phi$
using the frequency $\omega$ of the perturbing term  as a bifurcation parameter. 
Our results provide explanations of the frequency-locking phenomena reported by earlier authors, in particular by \cite{Rabinovich06}.

\begin{proposition}\label{prop:deltaLarge}
For every $\delta>\Phi$ and every $\gamma>0$ there is a value  $\omega_1>0$
such that for arbitrarily large $\omega> \omega_1$
there are pairs  of frequency locked periodic solutions of \eqref{general} appearing at saddle-nodes. Near the saddle-node  one of the pair of solutions is attracting.
Moreover, for $n\in\NN$
\begin{enumerate}
\item\label{item:k>=1}
if $k\ge 1$, then  the pairs of  solutions with frequency ratio $1:n$
occur for all  $\omega> n\omega_1$;
\item\label{item:k>1}
for $k>1$ when $\omega\to\infty$ each solution with frequency ratio $1:n$  
collapses into a constant solution;
\item\label{item:k=1}
for $k=1$ when $\omega\to\infty$ each pair of  solutions with frequency ratio $1:n$ 
collapses into a constant solution;
\item\label{item:k<1}
if  $0<k<1$, then  there is a value $\omega_2>\omega_1$ such that
 the pairs of solutions with frequency ratio $1:n$  
occur for  $n\omega_1<\omega<n\omega_2$.
\end{enumerate}
\end{proposition}
Note that the behaviour described in \eqref{item:k>1} can be seen as a Hopf bifurcation at $\omega=\infty$ (a periodic solution collapsing into a constant solution), 
while \eqref{item:k=1} would be a degenerate Hopf bifurcation (two periodic solutions collapsing into the same constant solution).

\begin{proof}
We consider separately the different values of $k$.

If $k\ge 1$, then it follows from Theorems~\ref{th:k>1} and \ref{th:k=1} that $\GG$ has two fixed points for each $\tau>\tau_1>0$ in regions \ref{item:C} and \ref{item:c}, respectively.
Proposition~\ref{prop:stabFold} states that one of these fixed points is attracting.
For $\omega> \omega_1=-K\pi/\ln\tau_1$ these fixed points correspond to frequency locked solutions with frequency ratio $1:1$
one of which is attracting.
Solutions with frequency ratio $1:n$
occur for $\omega>n \omega_1=-nK\pi/\ln\tau_1$ 
%and $\dpt \lim_{n\to\infty}\omega_{1,n}=\infty$ 
establishing assertion \eqref{item:k>=1}.

The two branches in region \ref{item:C} with $k>1$ terminate at $\tau=1$ at two different values of $s$.
This corresponds to two different  constant solutions of \eqref{general} when $\omega=-K\pi/\ln\tau$ tends to $\infty$ as in assertion \eqref{item:k>1}.
For assertion \eqref{item:k=1} note that when $k=1$ we are at region \ref{item:c} where there is a fold point at $\tau=1$ so the two periodic solutions collapse into the same constant solution.
\medbreak

Finally, concerning assertion \eqref{item:k<1}, if $0<k<1$
with  $\delta>\Phi$  we are in region \ref{item:Z} of Theorem~\ref{th:k<1}.
There are two fixed points of $\GG$ for each $\tau\in(\tau_1,\tau_2)$ with $0<\tau_1<\tau_2<1$.
Let $\omega_1=-K\pi/\ln\tau_1$ and $\omega_2=-K\pi/\ln\tau_2$.
Hence, for $\omega\in\left(n\omega_1,n\omega_2\right)$ 
these fixed points correspond to frequency locked solutions with frequency ratio $1:n$.
The subcritical bifurcation at $\tau_2$ includes an attracting fixed point by Proposition~\ref{prop:stabFold}.
\end{proof}

For $\delta<\Phi$
we do know that one solution branch bifurcating from each fold point is stable and then it loses stability at a Hopf bifurcation.
When 
$\delta>\Phi$  we have less information since away from the fold points there might be other Hopf bifurcations that we have not detected.

\section{Discussion and concluding remarks}
\label{s: discussion}
In this paper  we contribute to a theoretical understanding of the dynamics generated by 
 periodic perturbations of a  heteroclinic cycle in the attracting, clean and  robust heteroclinic network constructed by Aguiar \emph{et al} \cite{ACL06}. 
  We extend the work started in \cite{LR18,LR21, RL2014} to a more general setting.
 We are particularly interested in the  interplay between two  parameters $\gamma$ and $\delta$ describing, respectively, the amplitude of the perturbation and the saddle-value of the network that quantifies its attractiveness. 
 
 Using the analytical expression of the reduced first return map for trajectories near the  cycle
 deduced in \cite{LR18}  we have identified a sharp difference in the
  bifurcation diagrams and dynamical regimes arising when
  $\delta<\Phi$ and $\delta \geq \Phi$.
 This is one of the main novelties of the article.

The results have been interpreted back in terms of the original differential equation \eqref{general}
starting with a discussion of frequency locking and
periodic solutions with very long periods. 
We determined parameter values where these solutions arise from attracting branches of the bifurcation diagrams.
 These solutions would be easily seen in a numerical simulation.
 The invariant tori obtained here would be less easy to find numerically since they arise in a narrow region of parameter space and we have not determined their stability.
 This is true a fortiori for the chaotic regime obtained from the Bogdanov--Takens bifurcation as it takes place in an exponentially narrow parameter strip. 
 One would probably not find this numerically unless one knew it to be there in the first place.
 
 Invariant tori and chaotic dynamics for \eqref{general} were obtained in \cite{LR21} by different global methods.
For small $\omega>0$ persistent strange attractors were found in \cite{Rodrigues_SIAM}  under a condition involving the eigenvalues via the deformation of an invariant torus.
  These take place in regions of parameter space  different from those considered here.

When $k>1$ if the amplitude $\gamma$ of the forcing is small and if $\delta$ is above the threshold $\delta=\Phi$, then
 we  find a maximal  period $T>0$  below which 
 \eqref{general} has two $T$-periodic solutions of different stability. 
 These orbits arise from
 a resonant Arnold tongue.    
 Numerics for a similar system in \cite{DT3, TD1}  confirm  that the  first return map \eqref{eq:PeriodicSolutions} captures the dynamics quite well.
 
 It is possible that the two scenarios  (I) frequency-locked periodic orbits and (II) orbits that fluctuate close to the period expected from the time-averaged estimate \cite{LR23.1} 
overlap as $\omega$ is varied   giving rise to bistability as in \cite[Section 6.2]{LR18}.

Periodic perturbation of a similar network in the May--Leonard system \cite{ML75} is treated
 in  \cite{AHL2001,DT3,TD1} under the hypothesis  that the perturbation has very small amplitude.
  Afraimovich and collaborators in \cite{AHL2001} find,  for the first return map of a weakly attracting network, an invariant closed curve and, for perturbations with slightly larger amplitude,  dynamics conjugate to a shift on two symbols. 
 We obtain the same result  for different equations with similar dynamics by a different method:
 we locate points of Hopf and
 discrete-time Bogdanov--Takens bifurcation  in the case $\delta<\Phi$.
 
   Tsai and  Dawes \cite{TD1} give an asymptotic definition of weak and strong attraction for 
a cycle in a similar  network; here we have a precise threshold, the golden number $\Phi$ separating the two regimes.
  They also obtain numerically an invariant circle when  $\delta$ is close to 1.

The present work, together with \cite{LR18,LR21, LR23.1},
is part of a systematic study of the dynamics associated to a parametric periodically-forced vector field unfolding a heteroclinic cycle in an asymptotically stable heteroclinic network.
In \cite{LR21} the analysis in the extended phase space shows 
the existence of an invariant torus that is destroyed giving rise to rotational horseshoes and strange attractors. The original differential equation has autonomous and non-autonomous terms governed by independent parameters.

In the same context it is shown in \cite{LR23.1} that as the forcing frequency tends to infinity the dynamics reduces to that of a cycle in a
network under constant forcing, the constant being the average value of the forcing term. 
Under small constant forcing, each cycle in
the network breaks up into an attracting periodic solution that persists for periodic forcing of high frequency.
 This is the generic outcome of  a small autonomous perturbation of an attracting homo/heteroclinic cycle: a single attracting periodic orbit. 
 Our results show that this is not necessarily the case for a time-periodic perturbation.
 In this sense, the present article provides an analytic
  justification  for the findings of  \cite{AHL2001,DT3,TD1}.

\subsection*{Final remark}
The dynamics near the attracting network  $\Sigma_0$ has qualitative features in common with that near homoclinic orbits in general, at least in the regime with large $\delta$.

Many questions remain for future work; some are highlighted above. One obvious problem is  the computation of the Melnikov map associated to the   intersection of $W^u(\ww)$ and $W^s(\vv)$ in the same spirit of \cite[Section 7]{Mohapatra2015}. It would be of substantial interest to relate our bifurcations with those of the torus-breakdown and Arnold tongue bifurcations explored in \cite{TD1}. An important work on the topic would be to put all partial results  of the present paper and those of \cite{AHL2001, DT3, LR18, LR21, LR23.1,TD1} together.

 Generalising the results of the present article to equations where $\gamma(rt)$   is a generic time-varying input (with rate $r>0$) passing through a bifurcation point and forcing the network to loose stability  will be a major challenge. In the spirit of Ashwin \emph{et al} \cite{Ashwin2017}, $B$-tipping points might occur where the dynamics may change abruptly.
 We also would like to see  if these dynamical bifurcation diagrams and regimes (which individually are quite generic kinds of dynamical system) arise in other applied problems. 
A general discussion of dynamics arising in the context of bifurcations under periodic perturbation is presented in \cite{Simo2018} with additional numerical results in \cite{Simo2016}, where it is stated about global bifurcations  that the {\it ``challenge is to be able to predict for which parameter values they occur''.}

\section*{Acknowledgements}
The authors are grateful to an anonymous referee and the editor whose attentive reading and useful comments improved the final version of the article.

\end{document}